\documentclass[12pt]{article}

\usepackage{amstext    }
\usepackage{amsthm    }
\usepackage{a4}
\usepackage[mathscr]{eucal}
\usepackage{mathrsfs}

\usepackage{amsmath}
\usepackage{amssymb}
\usepackage{amscd}

\newtheorem{theorem}{Theorem}[section]
\newtheorem{definition}[theorem]{Definition}
\newtheorem{proposition}[theorem]{Proposition}
\newtheorem{corollary}[theorem]{Corollary}
\newtheorem{lemma}[theorem]{Lemma}
\newtheorem{remark}[theorem]{Remark}

\newtheorem{question}[theorem]{Question}

\newcommand{\cali}[1]{\mathscr{#1}}

\newcommand{\volume}{{\rm vol}}

\newcommand{\supp}{{\rm supp}}

\newcommand{\loc}{{loc}}
\newcommand{\ddc}{dd^c}
\newcommand{\dc}{d^c}
\newcommand{\dbar}{\overline\partial}
\newcommand{\ddbar}{\partial\overline\partial}

\newcommand{\PSH}{{\rm PSH}}
\newcommand{\PB}{{\rm PB}}

\newcommand{\id}{{\rm id}}
\newcommand{\pr}{{\rm pr}}

\newcommand{\Cc}{\cali{C}}
\newcommand{\Dc}{\cali{D}}

\newcommand{\Fc}{\cali{F}}

\newcommand{\Hc}{\cali{H}}
\newcommand{\Kc}{\cali{K}}

\newcommand{\Pc}{\cali{P}}
\newcommand{\Qc}{\cali{Q}}
\newcommand{\Rc}{\cali{R}}
\newcommand{\Uc}{\cali{U}}

\newcommand{\FS}{{\rm FS}}

\newcommand{\C}{\mathbb{C}}

\newcommand{\R}{\mathbb{R}}
\renewcommand{\P}{\mathbb{P}}

\title{Dynamics of horizontal-like maps in higher dimension}

\author{Tien-Cuong Dinh, Vi\^et-Anh Nguy{\^e}n and Nessim Sibony}

\begin{document}

\maketitle

\begin{center}
{\it Dedicated to Professor Gennadi Henkin\\ on the occasion of his 65th
  birthday}
\end{center}

\begin{abstract}
We study  the regularity of the Green currents and of the equilibrium measure 
associated to a horizontal-like map in $\C^k$, under
a natural assumption  on the  dynamical degrees. 
We estimate the speed of convergence  towards the Green currents, the decay of
correlations for the equilibrium measure
and the Lyapounov exponents. We show in particular that the
equilibrium measure is hyperbolic. We also show that the Green
currents are the unique invariant vertical and horizontal positive
closed currents.
The results apply, in
particular,  to H{\'e}non-like maps, to regular
polynomial automorphisms of $\C^k$ and to their small pertubations. 
\end{abstract}

\noindent
{\bf Classification AMS 2000:} Primary 37F, Secondary  32U40, 32H50.

\noindent
{\bf Keywords: } horizontal-like map, Green current, equilibrium
measure, entropy, mixing, Lyapounov exponent, structural disc of currents.

\section{Introduction} \label{introduction}

The abstract theory of non-uniformly hyperbolic systems is
well-developed, see e.g. Katok and   Hasselblatt \cite{kh},
Pesin \cite{pe},  L.-S. Young \cite{yo}.
It is  however difficult to show that a concrete  example is a non-uniformly 
hyperbolic system. The main  questions are to construct  
a measure of maximal entropy, to study  the decay of correlations
and to show  that the Lyapounov  exponents do  not vanish.  Such problems have been studied 
in  dimension $2$ for real H{\'e}non  maps by  Benedicks-Carleson, L.-S. Young, Viana etc.,
see e.g. \cite{bc, BonattiDiazViana, yo2}.
In this paper  we consider  these questions  for holomorphic horizontal-like maps  in $\C^k$ using  tools
from complex analysis: positive closed currents, estimates for solutions of
the $\ddc$-equation and appropriate spaces of test forms. The complex analytic methods permit to avoid
the delicate arguments used in the  real setting.

In \cite{ds1} the first  and the third authors  studied  the dynamics of polynomial-like maps
in several complex variables  using adapted spaces of test functions.
This approach permits to study  
convergence problems, in particular, the decay of correlations for the
measure of maximal entropy. 
Recall that a polynomial-like map is a proper holomorphic map
$f:U\rightarrow  V$ between convex open
sets  $U\Subset V$ (or more generally pseudoconvex open sets) in $\C^k$. Such  a map is somehow ``expanding",
but it has in general  a   non-empty critical set; so, it is not
uniformly hyperbolic  in the dynamical sense, see \cite{kh}. 
It is shown in \cite{ds1} that the measure of maximal entropy is
hyperbolic if the topological degree is strictly larger than  
the other dynamical degrees. This condition is natural and is stable under small
pertubations on the map. 
Holomorphic  endomorphisms of $\P^k$ can be lifted to polynomial-like maps in some open sets
of $\C^{k+1}$. So, their dynamical study is a special case of
polynomial-like maps. Small transcendental pertubations of such maps provide
large families of examples.

Here, we consider the quantitative aspects of the dynamics of horizontal-like maps $f$ in any
dimension, inside a product of convex open sets $D=M\times N$ in
$\C^p\times \C^{k-p}$. 
They are basically 
holomorphic maps  
which are somehow ``expanding" in $p$ directions (horizontal
directions) and 
``contracting" in the other $k-p$ directions (vertical directions), see Section  \ref{section_horizontal} for
the precise definition. They partially
look like a horseshoe. But, 
the expansion and contraction are of global nature, and in general, 
these maps are not uniformly hyperbolic. 
Small pertubations of horizontal-like maps are horizontal-like
provided that we shrink slightly the domain of definition.  
When $p=k$ we obtain polynomial-like maps.

H{\'e}non maps in $\C^2$ were studied by
Bedford-Lyubich-Smillie \cite{BedfordLyubichSmillie} with the
equilibrium measure introduced by the third author of the present
article, see also \cite{fs1}. The case of
horizontal-like maps in dimension 2, i.e. $k=2$ and $p=1$, has been  studied by Dujardin 
with emphasis on biholomorphic maps (H{\'e}non-like maps) \cite{duj}
and was developed by Dujardin, the first and the third authors to deal with random iteration of
meromorphic horizontal-like maps  \cite{dds}. It turns out that horizontal-like maps are  the building blocks 
for polynomial maps of ``saddle type". In particular, they were used 
to study rates of escape to infinity for polynomial mappings 
in $\C^2$. The randomness comes from the indeterminacy
points at infinity, see also \cite{Vigny}.

In this paper, we continue our study in the higher dimensional case.
In order to simplify the notation, {\bf we only consider invertible maps}.
However, a large part of our study can be extended to the general case.
Some basic objects and the first properties for such maps (Green currents
$T_\pm$, equilibrium measure $\mu$, entropy, mixing, 
etc) were constructed and established in \cite{ds6}.
The Green current $T_+$ is positive closed of bidegree $(p,p)$,
invariant under $f^*$ and is
{\it vertical}: its support does not intersect the vertical
boundary $\partial M\times N$ of $D$. The Green current $T_-$ is positive closed of bidegree $(k-p,k-p)$,
invariant under $f_*$ and 
is {\it horizontal}. The equilibrium measure $\mu$ is an invariant
probability measure which is equal to the wedge-product $T_+\wedge
T_-$ of the Green currents. The definition of wedge-product relies on
an intersection theory for positive closed currents. 

The main technical problem  is  the  use of currents
of bidegree $(p,p)$,  $p\geq 1$.  
For that purpose, 
a geometry on the space of positive closed $(p,p)$-currents 
was introduced using as basic objects: {\it structural discs of currents.}
Roughly speaking, in order to travel from  a positive 
closed current $R_1$ of bidimension $(k-p,k-p)$  to another one $R_2$, 
we construct a family of currents parametrized by a holomorphic disc $\Delta\subset\C$.
These currents appear as the slices of a  positive closed current
$\Rc$ of bidimension
$(k-p+1,k-p+1)$ in $\Delta\times D$; the currents $R_1$ and $R_2$ are seen as two points of
the disc, i.e. two currents obtained by slicing $\Rc$ with $\{\theta_1\}\times D$ and
$\{\theta_2\}\times D$ for some $\theta_1$, $\theta_2$ in $\Delta$.
We use properties of subharmonic functions on those
structural discs  in order to define the wedge-product of currents of
higher bidegree and
in order to prove the convergence results in the construction of
$T_\pm$ and $\mu$. More formally as in \cite{ds8} we use super-functions,
i.e. functions defined on horizontal currents which are p.s.h. on
structural discs of currents.

In the present article, we study the quantitative properties of these
basic dynamical objects. For a horizontal-like map $f$, one associates a
main dynamical degree $d\geq 2$ which is an integer. The topological
entropy of $f$ and the entropy of $\mu$ are equal to $\log d$. We will
define the other dynamical degrees $d_s^\pm$ in Section
\ref{section_horizontal}. 
One of our main results is the following.

\begin{theorem} \label{th_main}
Let $f$ be an invertible horizontal-like map on a convex domain
$D=M\times N$ in $\C^p\times \C^{k-p}$. Assume that the main dynamical
degree $d$ of $f$ is strictly larger than the other dynamical
degrees. Then the Green currents $T_+$ and $T_-$ of $f$ are the unique, up to
a multiplicative constant, invariant vertical and horizontal positive
closed currents of bidegrees $(p,p)$ and $(k-p,k-p)$ respectively.
The equilibrium measure $\mu$ of $f$ is exponentially
mixing and is hyperbolic. More precisely, $\mu$ admits $k-p$ strictly
negative  and $p$ strictly positive Lyapounov exponents.
\end{theorem}

We study the speed of convergence towards the Green currents
$T_\pm$ and
the equilibrium measure $\mu$, and also the regularity of these objects. 
The regularity is studied by considering on
which  space of forms  or functions  the currents or measures act
continuously.  We show in particular that
$\mu$  is $\PB,$ that is, plurisubharmonic functions
(p.s.h. for short)  are $\mu$-integrable.
The main tools here are estimates and localization of the support for good solutions
of the $\ddc$-equation. We obtain these estimates through integral
formulas (a classical result by
Andreotti-Grauert is crucial here). 
They permit to apply the $\ddc$-method 
and the duality method as in \cite{ds1,ds4,ds5,ds6}. The speed of
convergence towards the Green currents is a basic ingredient in the
proof of the decay of correlations for $\mu$. 

For H{\'e}non like-maps ($k=2$, $p=1$), the hypothesis on the dynamical
degrees is always satisfied. Theorem \ref{th_main}, except for the decay of
correlations (exponential mixing), was proved in \cite{duj}.
The decay of correlations  for H{\"o}lder observables and for H{\'e}non
maps was investigated by the first author in \cite{din}.
The hyperbolicity of the equilibrium measure is considered in a very
general context for meromorphic maps on compact K\"ahler manifolds by de
Th{\'e}lin \cite{deThelin}. We follow his method.

We end this introduction by giving another large family of examples.
Consider a polynomial automorphism $f$ of $\C^k$. 
We still denote by $f$ its meromorphic extension to $\P^k$.
When the indeterminacy sets $I_+$ and $I_-$ of $f$ and $f^{-1}$ 
in the hyperplane at infinity $L_\infty$ are non-empty 
and have no intersection,
we say that $f$ is {\it regular}. Then there is an integer $p$  such that
 $\dim I_{+}=k-p-1$ and $\dim I_{-}=p-1$. We refer to \cite{si} for
 the basic dynamical objects and properties of such maps, see also
 Section \ref{section_ex} below. 
Let $z=(z_1,\ldots,z_k)$ denote the coordinates in $\C^k$ and denote
$[z_0:\cdots:z_k]$ the homogeneous coordinates of $\P^k$. The hyperplane
 at infinity $L_\infty:=\P^k\setminus \C^k$ is given by the equation $z_0=0$.

\begin{corollary} \label{cor_ex}
Let $f$ be   a regular polynomial automorphism of $\C^k$. Assume
that the indeterminacy sets of $f$ and $f^{-1}$ are linear and defined by 
$$I_+=\{z_0=z_1=\cdots=z_p=0\} \quad \mbox{and}\quad 
I_-=\{z_0=z_{p+1}=\cdots=z_k=0\}.$$
Let $B_s^R$ denote the ball of center $0$ and of radius $R$ in $\C^s$.
Then, if $R$ is large enough, 
any holomorphic map on $B_p^R\times B_{k-p}^R$, close enough to 
$f$, is horizontal-like. Moreover, its equilibrium measure is exponentially
mixing and hyperbolic.
\end{corollary}

Note that the above pertubation of $f$ may be transcendental and that
Corollary \ref{cor_ex} produces large families of examples.

Here is a brief outline of the paper. In Section \ref{section_ddc},
the main tools, in particular, several classes of
currents and the solution of the $\ddc$-equation, are introduced. In Section \ref{section_horizontal}, we
recall the dynamical objects associated to a horizontal-like
map. Theorem \ref{th_main} is proved in Sections \ref{section_cv} and
\ref{section_measure}. Corollary \ref{cor_ex} is deduced from Theorem \ref{th_main} and 
from Proposition \ref{prop_ex} in the last section. Also in the last
section open questions are
stated.

\medskip
\noindent
{\bf Notation and convention.} Throughout the paper, $D:=M\times N$ is a bounded convex domain
in $\C^p\times \C^{k-p}$. 
The estimates we obtain are valid in the
interior of $D$ and might be bad near the
boundary, but this is harmless for the
type of maps we consider. So, we sometimes reduce $D$ slightly 
in order to have maps and currents defined in a neighbourhood of
$\overline D$; this simplifies the exposition.  We will also choose strictly convex domains
with smooth boundary $M''\Subset M'\Subset M$ and $N''\Subset N'\Subset
N$ and consider the domains 
$D':=M'\times N'$ and $D'':=M''\times N''$.
When we consider
vertical currents $R$ or horizontal currents $S$, $\Phi$, our choice
is so that $R$ is supported on $M''\times N$ and
$S$, $\Phi$ are supported on $M\times N''$. When we consider a
horizontal-like map $f$ on $D$, we assume that $f^{-1}(D)\subset
M''\times N$ and $f(D)\subset M\times N''$. So, $f$ restricted to $D'$
or $D''$ is horizontal-like. The convex domains $\widetilde M$,
$\widehat M$, $\widetilde N$, $\widehat N$ are chosen so that
$M\Subset \widetilde M\Subset \widehat M$ and $N\Subset \widetilde
N\Subset \widehat N$. Note also that {\it when we consider the convergence of a family of
vertical or horizontal currents, we assume that they have support in the same
vertical or horizontal set.}


\section{Currents and $\ddc$-equation} \label{section_ddc}

In this section, we will introduce the tools used in this work. We
will give some geometrical and analytical properties of several
classes of currents. In particular, we will define structural discs of
currents and solve the $\ddc$-equation with estimates and with controlled support.  
Recall that $\dc:={i\over 2\pi} (\dbar -\partial)$.

\medskip

\noindent
$\bullet$ {\bf Vertical, horizontal currents and their intersection.}
We call {\it vertical} 
(resp. {\it horizontal}) 
{\it boundary} of $D$ the sets $\partial_v D:=\partial M\times N$ 
(resp. $\partial_h D:=M\times \partial N$). 
A subset $E$ of $D$ is  {\it vertical} (resp. {\it horizontal}) if
$\overline E$ does not intersect  $\overline{\partial_v D}$ 
(resp. $\overline{\partial_h D}$). Let $\pi_1$ and $\pi_2$ denote the canonical 
projections of $D$ onto $M$ and $N$. Then $E$ is vertical or horizontal if and only if 
$\pi_1(E)\Subset M$ or $\pi_2(E)\Subset N$.
A current on $D$ is {\it vertical} or
{\it horizontal} if its support is vertical or horizontal. 
Let $\Cc_v(D)$ denote 
the cone of positive closed vertical
currents of bidegree $(p,p)$ on $D$.
Consider a current $R$ in $\Cc_v(D)$. 
Since $\pi_2$ is proper on $\supp(R)$, $(\pi_2)_*(R)$ is a positive closed current of
bidegree $(0,0)$ on $N$. Hence, $(\pi_2)_*(R)$ is given by a constant function on
$N$ that we denote by  $\| R\|_v$. Convergence in $\Cc_v(D)$ is the
weak convergence of currents with support in a fixed vertical set.

Recall from Theorem 2.1 in \cite{ds6} that  the slice measure  $\langle R,\pi_2,w\rangle$
 is defined for every $w\in N$, and that 
its mass  is equal to  $\| R\|_v$ which is independent of $w$. 
We say that  $\| R\|_v$ is the {\it slice mass} of $R$.
For every smooth probability measure $\Omega$ with compact support in $N$, we have 
$\| R\|_v:=\langle R,(\pi_2)^*(\Omega)\rangle$.  
When $\| R\|_v=1$ we say that $R$ is {\it normalized}. Let $\Cc_v^1(D)$ denote the set 
of such currents. This convex set is relatively compact in the cone of
positive closed currents on $D$. In particular, the mass of 
normalized currents $R$ on  a compact set of $D$ is bounded uniformly
on $R$. In order to avoid convergence problems on the boundary, we will also
use the convex set $\Cc_v^1(M\times \overline N)$ of positive closed currents
which are vertical in $M\times \widetilde N$ with slice mass 1 for some neighbourhood
$\widetilde N$ of $\overline N$.

The slice mass $\|\cdot\|_h$, the sets $\Cc_h(D)$, $\Cc_h^1(D)$ and
the convergence for 
horizontal currents of bidegree $(k-p,k-p)$ are  defined similarly.
If $R$ is a current in $\Cc_v(D)$ and $S$ is a current in $\Cc_h(D)$
we can define the intersection $R\wedge S$. This is a positive measure
of mass $\|R\|_v\|S\|_h$ with support in $\supp(R)\cap\supp(S)$, see
\cite{ds6}. It depends linearly on $R$ and on $S$ and is
continuous with respect to the {\it plurifine topology} in the following sense. Let
$(R_\theta)$ and $(S_{\theta'})$ be structural discs in $\Cc^1_v(D)$
and $\Cc^1_h(D)$, see the definition below. Assume that
$\supp(R_\theta)\cap\supp(S_{\theta'})$ is  contained in an open set
$\Omega\Subset D$. If $\varphi$ is a p.s.h. function on a neighbourhood
of $\overline \Omega$, then $\langle R_\theta\wedge
S_{\theta'},\varphi\rangle$ is either a p.s.h. function of
$(\theta,\theta')$ or equal identically to $-\infty$, see Proposition 3.4 and Remark 3.8 in \cite{ds6}.
Basically, for a suitable choice, with $R_1=R$, $S_1=S$ and $R_\theta$, $S_{\theta'}$
smooth when $\theta\not=1$, $\theta'\not=1$, we
obtain $R\wedge S$ as the limit of $R_\theta\wedge S_{\theta'}$,
$R_\theta\wedge S$, $R\wedge S_{\theta'}$ for
$\theta\rightarrow 1$ and $\theta'\rightarrow 1$.  
It is also shown in \cite{ds6} that for a p.s.h. function $\varphi$ on $D$
$$\langle R\wedge S,\varphi\rangle
=\limsup \langle R'\wedge S',\varphi\rangle
=\limsup \langle R'\wedge S,\varphi\rangle
=\limsup \langle R\wedge S',\varphi\rangle,$$
with $R'$, $S'$ smooth in $\Cc_v(D)$, $\Cc_h(D)$ converging
respectively to $R$ and $S$.

\medskip

\noindent
$\bullet$ {\bf Structural discs of currents.} Let $X$ be a complex manifold.
Consider a positive closed $(p,p)$-current $\Rc$ in $X\times D$. We
assume that the support of $\Rc$ is contained in $X\times M'\times
N$ for some open set $M'\Subset M$. Let $\pi:X\times
D\rightarrow X$ denote the canonical projection. It is shown in
\cite{ds6} that the slice $\langle\Rc,\pi, x\rangle$ exists for
every $x\in X$. They can be considered as the intersection of
$\Rc$ with the current of integration on $\pi^{-1}(x)$. This is a
positive closed $(p,p)$-current on $\{x\}\times D$ that we
identify with a current $R_x$ on $D$ which is vertical. 
When $\Rc$ is a smooth form, the slice $R_x$ is simply the restriction
of $\Rc$ to $\pi^{-1}(x)$.
The slice mass of $R_x$ does not depend on $x$. So,
multiplying $\Rc$ with a constant, we can assume that this mass is 1.
We obtain a map $\tau:X\rightarrow\Cc^1_v(D)$ with
$\tau(x):= R_x$. In general, $R_x$ does not depend
continuously on $x$ with respect to the usual topology on $X$. The dependence is
continuous with respect to the {\it plurifine topology}, i.e. the coarsest topology for which
p.s.h. functions on $X$ are continuous. We call 
{\it structural variety} of $\Cc_v^1(D)$ the map $\tau$ or
the family $(R_x)$. This notion can be easily extended to
$\Cc_v^1(M\times\overline N)$.

Consider a vertical positive closed $(p,p)$-current $R$ in
$\Cc_v^1(M\times\overline N)$. So, $R$ is a vertical current of
slice mass 1 on $M'\times \widetilde N$ for some convex open sets
$M'\Subset M$ and 
$\widetilde N\Supset N$. 
Let $\Delta$ denote a small neighbourhood of the interval $[0,1]$ in $\C$. 
We constructed in \cite{ds6} a {\it particular structural
disc} $(R_\theta)_{\theta\in\Delta}$ in $\Cc_v^1(M\times \overline N)$
parametrized by $\Delta$ such that $R_1=R$ and $R_0$ is independent of
$R$.
The current $R_\theta$ is obtained as a regularization of $R$. More
precisely, we consider some holomorphic family of linear endomorphisms $h_{a,b,\theta}:\C^k\rightarrow\C^k$
parametrized by $(a,b,\theta)\in \C^p\times\C^{k-p}\times\Delta$ with $h_{a,b,1}=\id$. The current $R_\theta$ is obtained
using a smooth probability measure $\nu$ with compact support in 
$\C^p\times\C^{k-p}$:
$$R_\theta:=\int (h_{a,b,\theta})_*(R)d\nu(a,b).$$
The convexity of $M\times N$ and the fact that $R$ is defined on $M'\times \widetilde N$ permit
to define the smoothing and to obtain vertical currents
$R_\theta$ in $\Cc_v^1(M\times \overline N)$. The size of $\Delta$ depends only on $M$, $M'$, $N$
and $\widetilde N$.
The considered structural discs 
satisfy the following important properties. {\it The currents $R_\theta$
depend continuously on $\theta$, linearly on $R$ and are smooth for
$\theta\not=1$. The continuity is with respect to the weak topology on
$R_\theta$ and the usual topology on $\theta$. 
Moreover,
 $R_\theta$ depend continuously on $\theta$ and on $R$ with respect to the usual topology on
$\theta\in\Delta\setminus\{1\}$, 
the $\Cc^\infty$ topology on $R_\theta$ and the weak topology on
$R$. When $R$ is smooth, the last property also holds for 
$\theta\in\Delta$.}

\medskip

\noindent
$\bullet$ {\bf PSH currents and p.s.h. functions.}
A real $(k-p,k-p)$-horizontal current  $\Phi$  on $D$  is 
called {\it PSH} if $\ddc\Phi\geq 0$\footnote{In other situations, we often assume that $\Phi$ is of order 0
or negative. This is necessary in particular when one defines the
pull-back by a non-invertible map \cite{ds7}. Note that a
p.s.h. function is defined everywhere but not a PSH current.}. 
Let $\PSH_h(D)$ denote the set of horizontal PSH 
currents. It is endowed with the following topology. A sequence
$(\Phi_n)$ converges to $\Phi$ in $\PSH_h(D)$ if $\Phi_n\rightarrow \Phi$
weakly and if $\Phi_n$ and $\Phi$ have their supports in a fixed horizontal
set of $D$.

Recall that an upper semi-continuous function $\phi:D\rightarrow
\R\cup\{-\infty\}$ is p.s.h. if it is not identically $-\infty$ and if
its restriction to any holomorphic disc in $D$ is subharmonic or
equal to $-\infty$. Let
$\PSH(D)$ denote the cone of such functions. It is relatively compact
in $L_\loc^p(D)$ for $1\leq p<+\infty$. Note that an $L^1_\loc$ function 
 $\phi:D\rightarrow
\R\cup\{-\infty\}$ is p.s.h. if it is strongly upper semi-continuous and if 
$\ddc\phi$ is a positive closed current. The strong upper semi-continuity means
$\phi(a)=\limsup_{z\rightarrow a}\phi(z)$ for $a\in D$ and $z\in A$
where $A$ is any measurable subset of full measure in $D$.
Denote by $\PSH(\overline D)$
the cone of p.s.h. functions defined in a neighbourhood of $\overline D$.

\medskip

\noindent
$\bullet$ {\bf Extension of spaces of test forms and super-functions.} Let $R$ be a
current in $\Cc_v^1(D)$. It acts on horizontal 
smooth forms of bidegree $(k-p,k-p)$. We will extend this space of
test forms. Let $\Hc_h(D)$ denote the space of 
real horizontal currents $\Phi$ of bidegree $(k-p,k-p)$ with $\ddc
\Phi=0$. We consider the following topology on $\Hc_h(D)$: a sequence
$(\Phi_n)$ converges to $\Phi$ in $\Hc_h(D)$ if
$\Phi_n\rightarrow\Phi$ weakly and  $\Phi_n$ have
support in a fixed horizontal set.

\begin{proposition} \label{prop_test_closed}
The action of $R$ can be extended in a unique way to a positive continuous
linear form on $\Hc_h(D)$. Moreover, $(R,\Phi)\mapsto \langle
R,\Phi\rangle$ with $\Phi\in\Hc_h(D)$ is bilinear and
continuous in $(R,\Phi)$. In particular, $\langle R,\Phi\rangle$ is
bounded on compact subsets of $\Cc^1_v(D)\times \Hc_h(D)$. 
\end{proposition}
\proof
Observe that if $\Phi$ is a current in $\Hc_h(D)$ we can use a slight
dilation and a convolution in order to regularize $\Phi$. 
So, there are smooth forms $\Phi_n$ converging to $\Phi$ in $\Hc_h(D)$.
This implies the uniqueness, the linearity and the positivity of the
extension. Recall that the positivity means $\langle R,\Phi\rangle\geq
0$ for $\Phi\geq 0$. We prove now the existence of the 
extension on $\Hc_h(D)$ and the continuity.

Shrinking $D$ allows to assume that $R$ is defined on $M'\times
\widetilde N$ with $M'\Subset M$ and $\widetilde N\Supset N$.
Consider the structural disc $(R_\theta)$ as above. Define 
$h(\theta):=\langle R_\theta,\Phi\rangle$. As in \cite[Thm. 2.1]{ds6},
$h$ is a harmonic function on $\Delta\setminus\{1\}$. If $\Phi$ is
smooth, the function is defined and is harmonic on $\Delta$. 
Define $h_n(\theta):=\langle R_\theta,\Phi_n\rangle$. The above
description of properties of $R_\theta$ implies that $h_n$ converge
locally uniformly to $h$ on $\Delta\setminus\{1\}$. Since
$h_n$ are harmonic on $\Delta$ and locally uniformly bounded, by
maximum principle, the limit $h$ can be
extended to a harmonic function on $\Delta$ and $h_n$ converge to $h$
on $\Delta$. Observe that the limit does not depend on the choice of
$\Phi_n$. 

We have $\langle R,\Phi\rangle =h(1)$ when $\Phi$ is smooth.
Define $\langle R,\Phi\rangle:=h(1)$ the extension of
$R$ to all $\Phi$ in $\Hc_h(D)$. 
Recall that $R_\theta$, for $\theta\not=1$, depends continuously on $R$ with
respect to the $\Cc^\infty$
topology on $R_\theta$. Hence, $h$ depends continuously on $(R,\Phi)$.
The continuity of  $\langle R,\Phi\rangle$ follows.
\endproof

We will extend $R$ to a linear form on $\PSH_h(D)$,
but the
extension can take the value $-\infty$. Recall that $R$ is a current
on $M'\times\widetilde N$.

\begin{proposition} \label{prop_test_psh}
The limit  
$\langle R,\Phi\rangle:=\limsup \langle R,\Phi'\rangle$ with $\Phi'$
smooth converging to $\Phi$ in $\PSH_h(M'\times N)$,
defines an extension of $R$ to  $\PSH_h(D)$.
The extension depends linearly on
$R$, $\Phi$. It takes values in $\R\cup\{-\infty\}$ and 
does not depend on the choice of $M'$ and $\widetilde N$. The function
$\theta\mapsto \langle R_\theta,\Phi\rangle$ is subharmonic on
$\Delta$ and we have
$\langle R,\Phi\rangle=\limsup \langle R',\Phi\rangle$
with $R'\rightarrow R$ in $\Cc_v^1(M\times \overline N)$. 
\end{proposition}
\proof
We can 
assume that $\Phi$ is supported on $\widetilde M\times N'$ and that
$R$ is vertical in $M''\times\widetilde N$. So, we can assume that the
considered currents $\Phi'$ are horizontal on $D$.
Consider first the case where $\Phi$ is smooth. Let 
$\Phi_n$ be a sequence of smooth
forms converging to $\Phi$ in $\PSH_h(D)$. Define $h(\theta):=\langle
R_\theta,\Phi\rangle$ and  $h_n(\theta):=\langle
R_\theta,\Phi_n\rangle$. These functions are subharmonic and
continuous on $\Delta$, see \cite[Thm 2.1]{ds6} (the subharmonicity
is deduced from the positivity of 
$\ddc (\Rc\wedge \Phi_n)$ and of its push-forward to $\Delta$). We also have $h_n\rightarrow h$ on
$\Delta\setminus\{1\}$. It follows from the classical Hartogs' lemma
\cite{Hormander} that
$\limsup h_n(1) \leq h(1)$. So, 
$\limsup \langle R,\Phi'\rangle\leq \langle R,\Phi\rangle$. 

On the other hand, since $R_\theta$ is obtained from $R$ by 
smoothing using an averaging on a group of linear transformations, a coordinate change implies that $\langle
R_\theta,\Phi\rangle =\langle R,\Phi_\theta\rangle$ where
$\Phi_\theta$ is obtained from $\Phi$ by a similar smoothing. The fact that
$\Phi$ is defined on $\widetilde M\times N'$ guarantees that
$\Phi_\theta$ is horizontal in $D$. We also have
$\Phi_\theta\rightarrow \Phi$ when $\theta\rightarrow 1$ for the
$\Cc^\infty$ topology.
Since $h$
is continuous we deduce that $\langle R,\Phi_\theta\rangle\rightarrow \langle
R,\Phi\rangle$ when $\theta\rightarrow 1$. So, $\langle R,\Phi\rangle= \limsup
\langle R,\Phi'\rangle$ when $\Phi$ is smooth. 
In other words, $\langle R,\Phi\rangle:=\limsup \langle R,\Phi'\rangle$
defines an extension of $R$ to all $\Phi$ in $\PSH_h(D)$. It is clear
that the extension does not depend on the choice of $\widetilde N$.

For a general current $\Phi$, there are smooth forms $\Phi_n$
converging to $\Phi$ in $\PSH_h(D)$. Define $h_n$ and $h$ as above. The
function $h$ is defined on $\Delta\setminus \{1\}$. The functions $h_n$
are continuous subharmonic, bounded from above and converge to $h$ on
$\Delta\setminus\{1\}$. It follows that $h$ can be extended to a
subharmonic function on $\Delta$. By Hartogs' lemma, we have
$$h(1)\geq \limsup h_n(1)=\limsup \langle R,\Phi_n\rangle.$$
It follows that $h(1)\geq \langle R,\Phi\rangle = \limsup
\langle R,\Phi'\rangle$. 

On the other hand, since $h$ is subharmonic, 
we have $h(1)=\limsup h(\theta)=\limsup\langle R,\Phi_\theta\rangle$ 
when $\theta\rightarrow 1$. We deduce
as above that $h(1)= \langle R,\Phi\rangle$. Since $h$ depends
linearly on $R$ and $\Phi$, $\langle R,\Phi\rangle$ depends linearly
on $R$ and $\Phi$. We also obtain that $\theta\mapsto \langle
R_\theta,\Phi\rangle$ is subharmonic on $\Delta$.

It remains to prove that $\langle R,\Phi\rangle=\limsup \langle R',\Phi\rangle$
with $R'\rightarrow R$ in $\Cc_v^1(M\times\overline N)$. 
This property implies that $\langle R,\Phi\rangle$ is independent of
the choice of $M'$. 
Since $\langle
R,\Phi\rangle=\limsup \langle R_\theta,\Phi\rangle$ for
$\theta\rightarrow 1$, we have $\langle R,\Phi\rangle\leq \limsup \langle R',\Phi\rangle$
with $R'\rightarrow R$. Now, if $(R'_\theta)$ is the structural disc
associated to $R'$ and if $h'(\theta):=\langle R'_\theta,\Phi\rangle$,
then $h'(\theta)\rightarrow h(\theta)$ for $\theta\not=1$. We deduce
from Hartogs' lemma that $h(1)\geq \limsup h'(1)$ which implies that 
$\langle R,\Phi\rangle\geq \limsup \langle R',\Phi\rangle$ and
completes the proof.
\endproof

\begin{remark} \label{rk_integral}
\rm
We can consider $R$ as a vertical current and $\Phi$ as a horizontal
one in appropriate domains $D'\Subset D$
and define $\langle R,\Phi\rangle$ on $D'$ instead of $D$. We will obtain the same
value. Indeed, 
in order to define
$(R_\theta)$ we can find smoothings which are adapted for
both $D$ and $D'$, see \cite{ds6} for details. 
\end{remark}

\begin{remark} \label{rk_int_psh}
\rm
Let $R$ be a current in $\Cc_v(D)$, $S$ in $\Cc_h(D)$ and $\varphi$ a
p.s.h. function on $D$. If $\varphi$ is integrable with respect to the
trace measure $S\wedge\omega^p$ of $S$ then $\varphi S$ defines a
current in $\PSH_h(D)$. We deduce from the above results that
$$\langle R\wedge S,\varphi\rangle=\limsup_{\theta\rightarrow 1}
\langle R_\theta\wedge S,\varphi\rangle=\limsup_{\theta\rightarrow 1}
\langle R_\theta,\varphi S\rangle=\langle R,\varphi S\rangle.$$ 
\end{remark}

\begin{definition} \rm
Let $\Lambda:\Cc^1_v(M\times \overline N)\rightarrow \R\cup\{-\infty\}$ be an upper
semi-continuous function
which is not identically $-\infty$. We say that $\Lambda$ is
a {\it p.s.h. super-function} if it is p.s.h. or identically equal to $-\infty$ on
each structural variety in $\Cc_v^1(M\times \overline N)$, and $\Lambda$ is {\it
  pluriharmonic} if both $\Lambda$ and $-\Lambda$ are p.s.h.,
see also \cite{ds8}.
\end{definition}

\begin{proposition} \label{prop_super_function}
Let $\Phi$ be a real horizontal $(k-p,k-p)$-current
  on $D$.  
If $\Phi$ is $\ddc$-closed, then $R\mapsto \langle
R,\Phi\rangle$ defines a pluriharmonic super-function.
If $\Phi$ is PSH, then $R\mapsto \langle
R,\Phi\rangle$ is a p.s.h. super-function.
\end{proposition}
\proof
We only have to prove the second assertion. 
Consider a structural variety $(R_x)_{x\in X}$ as above. 
Without loss of generality, we can assume that $R_x$ are vertical in
$M'\times \widetilde N$ and $\Phi$ is horizontal in $M\times N'$, see
Remark \ref{rk_integral}. 
We want to
prove that $x\mapsto \langle R_x,\Phi\rangle$ is identically equal to
$-\infty$ or p.s.h. If $\Phi$ is smooth, this
was proved in \cite[Lemma 2.2]{ds6}. 
For the general case, we have
$$\langle R_x,\Phi\rangle =\limsup_{\theta\rightarrow 1} \langle
R_{x,\theta},\Phi\rangle =\limsup_{\theta\rightarrow 1} \langle
R_x,\Phi_\theta\rangle,$$
where $(R_{x,\theta})$ is the particular structural disc constructed
as above using the same smoothing for each $R_x$.  
We deduce from the regularity of $R_{x,\theta}$ that $\langle
R_{x,\theta},\Phi\rangle$ is locally uniformly bounded on
$(x,\theta)\in X\times (\Delta\setminus\{1\})$. Since $\theta\mapsto \langle
R_x,\Phi_\theta\rangle$ is p.s.h., 
it follows from the
maximum principle that $\langle R_{x,\theta},\Phi\rangle$ is locally uniformly
bounded from above on $X\times\Delta$. Hence, the upper
semi-continuous regularization of $x\mapsto \langle R_x,\Phi\rangle$
is p.s.h. or identically $-\infty$. 
It is enough to show that 
$x\mapsto \langle R_x,\Phi\rangle$ is 
upper semi-continuous. 

For every $a\in X$, we have $\limsup_{x\rightarrow a}\langle
R_x,\Phi_\theta\rangle =\langle R_a,\Phi_\theta\rangle$ for
$\theta\not=1$. Since the functions $\theta\mapsto \langle R_x, \Phi_\theta\rangle$
are subharmonic, we deduce using Hartogs' lemma that $\limsup_{x\rightarrow a}\langle
R_x,\Phi_\theta\rangle \leq \langle R_a,\Phi_\theta\rangle$ for every
$\theta$. This implies the result. 
\endproof

\medskip

\noindent
$\bullet$ {\bf PB, PC currents and measures.} 
Let $T$ be a vertical current of bidegree $(p,p)$ in $\Cc_v(D)$.
We say that $T$ is {\it PB} if $\langle T,\Phi\rangle$ is bounded when
$\Phi$ is in a relatively compact subset of $\PSH_h(D)$.
We say that $T$ is {\it PC}
if it can be extended to a continuous linear form on $\PSH_h(D)$ with
respect to the topology we have introduced.
Observe that this extension coincides with the extension in Proposition
\ref{prop_test_psh}.
PC currents are PB.
PB and PC horizontal currents of bidegree $(k-p,k-p)$ are defined in
the same way. In the case of bidegree $(1,1)$, PB and PC currents
correspond to currents with bounded and continuous local potentials,
see also \cite{ds1, ds4, ds5}.

A positive measure $\mu$ with compact support in $D$ is said to be
{\it PB} if  $\langle \mu,
\phi\rangle $ is bounded when $\phi$ are smooth functions in a
relatively compact subset in $\PSH(\overline D)$.
Since p.s.h. functions on a neighbourhood of $\overline D$ can be
approximated by decreasing sequences of smooth ones, $\mu$ is PB if and only if 
p.s.h. functions on a neighbourhood of $\overline D$ are $\mu$-integrable.
PB measures have no mass on {\it pluripolar sets}, i.e. sets
which are contained in the pole set $\{\phi=-\infty\}$ of a
p.s.h. function $\phi$.
The measure $\mu$ is said to be {\it PC} if it can be extended to a linear
continuous form on $\PSH(D)$. Denote by
$\langle \mu,\phi\rangle$ the value of this extension on $\phi$. Note
that by continuity the extension is unique and
$\langle\mu,\phi\rangle$ is equal to the usual 
integral $\langle \mu,\phi\rangle$ of 
$\phi$. Any PC measure is PB.

\medskip

\noindent
$\bullet$ {\bf Solution of $\ddc$-equation.} We consider the
$\ddc$-equation on $D$. We will need negative solutions with
horizontal or vertical support and with estimates on the mass. The
behavior near the rest of the boundary is not important in our
study. 
The following theorem is obtained using
classical results. Recall that $\dc:={i\over 2\pi}(\dbar-\partial)$. 

\begin{theorem} \label{th_ddc}
Let $M'$ and $M$ be convex domains 
in $\C^p$ such that $M'\Subset  M$. Let $N'$ and $N''$ be convex open sets in 
 $\C^{k-p}$ such that $N''\Subset N'$. Let $\Omega$ be a
 horizontal positive closed current of bidegree $(k-p+1,k-p+1)$ on
 $M\times N''$. Then there is a horizontal negative $L^1$
 form $\Phi$ of bidegree $(k-p,k-p)$ on $M'\times N'$ such
 that 
$$\ddc \Phi=\Omega \quad \mbox{on} \quad   M'\times N'\qquad
\mbox{and} \qquad
\|\Phi\|_{M'\times N'}\leq c\|\Omega\|_{M\times  N''}$$ 
with $c>0$ independent of $\Omega$.
Moreover, $\Phi$ is defined by an integral formula, and depends
linearly and continuously on $\Omega$. 
\end{theorem}

In what follows, the solutions of $d$, $\dbar$ or $\ddc$ equations
are given by classical integral formulas. Consequently, the linearity, the continuous
dependence on data and the estimate on the mass of solutions are 
satisfied. Therefore, we will focus our attention only on the support of the solutions.

\begin{lemma}
Let $D'$ and $D$ be convex domains in $\C^k$ with $D'\Subset
D$. Let $\Omega$ be a positive closed current of bidegree
$(k-p+1,k-p+1)$ on $D$. There is a  negative $L^1$
 form $\Psi$ of bidegree $(k-p,k-p)$ on $D'$, smooth out of the
 support of $\Omega$, such
 that $\ddc \Psi=\Omega$ on $D'$.
\end{lemma}
\proof
We can assume that $D$ is contained in the ball of
center 0 and of radius $1/2$.
Define for coordinates $(z,\xi)$ on $\C^k\times \C^k$ the kernel
$$K(z,\xi):=\log\|z-\xi\|(\ddc \log\|z-\xi\|)^{k-1}.$$
Observe that $K$ is negative when $\|z\|<1/2$, $\|\xi\|<1/2$, and $\ddc K$ is equal to the current of integration on the
diagonal of $\C^k\times\C^k$. Let $\chi$ be a cut-off function, $0\leq
\chi\leq 1$, with compact support in $D$,
such that $\chi=1$ on a neighbourhood $U$ of $\overline D'$. Define 
$$\Psi'(z):=\int_\xi \chi(\xi)\Omega(\xi)\wedge K(z,\xi).$$
Hence, $\Psi'$ is a negative $L^1$ form depending continuously on $\Omega$. If $z$ is outside the support of
$\Omega$, then $\Psi'(z)$ is given by an integration outside the
singularities of $K$. So, $\Psi'(z)$ is smooth there.

Let $\pi_1$ and $\pi_2$ denote the canonical projections of
$\C^k\times\C^k$ on its factors. If $\Omega$ is smooth we have
$$\Psi'=(\pi_1)_*(\pi_2^*(\chi\Omega)\wedge K).$$
Since $\Omega$ is closed and $\ddc K=[z=\xi]$, we deduce that
$\Omega':=\ddc\Psi'-\Omega$ is equal on $U$ to 
\begin{eqnarray*}
\Omega' & = & \int_\xi d\chi(\xi)\wedge \Omega(\xi)\wedge\dc K(z,\xi) -\int_\xi
\dc\chi(\xi)\wedge\Omega(\xi)\wedge dK(z,\xi)\\
& & +\int_\xi\ddc\chi(\xi)\wedge\Omega(\xi)\wedge K(z,\xi).
\end{eqnarray*}
The last formula is valid for arbitrary $\Omega$ by regularization. 
So, $\Omega'$ is defined by integration on
$\{d\chi(\xi)\not=0\}$ where $K(z,\xi)$ is smooth if $z\in U$. It
follows that $\Omega'$ is smooth. We also have good
estimates on $\Cc^r$ norm of this form on compact subsets of $U$.

Since $\Omega'$ is closed and smooth, it is classical to obtain smooth solution of the equation
$\ddc \Psi''=\Omega'$ with estimates (we first solve 
a $d$-equation and then a $\dbar$-equation, the method will be described
below with details in a situation where more estimates are needed). One checks
 that $\ddc\Psi=\Omega$ for
$\Psi:=\Psi'-\Psi''-c\omega^{k-p}$ where $\omega:=\ddc\|z\|^2$ is the
standard K{\"a}hler form on $\C^k$ and $c>0$ is large enough in order to
guarantee that $\Psi$ is negative on $D'$.
\endproof

Now, we  need to control the support of the solution.
We shrink slightly $M$ and extend slightly $N''$. This allows
to assume that $\Omega$ is defined in $\widetilde M\times F$ for some fixed
compact set $F$ in $N''$.
Using the previous lemma, we can find $\Psi$ on $M\times
N$,  smooth outside the support of $\Omega$
such that $\ddc\Psi=\Omega$. Let $\chi$ be a
cut-off function equal to 1 on a neighbourhood of
$M\times F$ and equal to 0 near
$M\times \partial N''$ and on $M\times (N\setminus N'')$. In particular, $\chi=1$ on the
support of $\Omega$ and $\Psi$ is smooth on $\{d\chi\not=0\}$. 
Define $\Phi_1:=\chi\Psi$ and $\Omega':=\ddc\Phi_1-\Omega$. This is a
smooth horizontal closed form of bidegree $(k-p+1,k-p+1)$ with support in
$M\times N''$. Moreover, $\Omega'$ vanishes near $M\times F$ and has a
controlled  $\Cc^r$ norm.
We will find a smooth positive solution of the equation $\ddc
\Phi_2=\Omega'$ with horizontal support in $M'\times N'$. The current
$\Phi:=\Phi_1-\Phi_2$ satisfies Theorem \ref{th_ddc}.

A construction using an integral formula as in the book
\cite[pp. 37-39 and 61-63]{bt} by Bott
and Tu implies that there is a real smooth form $\Psi$ which is horizontal
in $M\times N''$ such that $d\Psi=\Omega'$ (shrink
$M$ and extend $N''$ if necessary). Of course, it satisfies
the desired estimates in $\Cc^r$ norms.
Moreover, we can write $\Psi=\Psi'+\Psi''$ with $\Psi'$ of bidegree
$(k-p,k-p+1)$ and $\Psi''$ of bidegree $(k-p+1,k-p)$ such that 
$\Psi''=\overline \Psi'$. 

\begin{lemma}
There is a smooth horizontal form $\Phi'$ on $M'\times N'$, 
of bidegree $(k-p,k-p)$,
such that $\dbar \Phi'=\Psi'$. 
\end{lemma}
\proof
Recall that we can, in each step of the proof, shrink or extend slightly the considered
domains $M$, $N'$ or $N''$. This permits to avoid the
problem near the boundary and to assume that they are strictly convex
with smooth boundary.
Since $d\Psi$ is of bidegree $(k-p,k-p)$, we have $\dbar\Psi'=0$. So,
using a classical integral formula (see, for example  \cite{hk1, ru}) we can find
a smooth form $\Phi^*$ of bidegree $(k-p,k-p)$ on $M\times
N$ such that $\dbar \Phi^*=\Psi'$. Its support is not necessarily
horizontal. So, we have $\dbar\Phi^*=0$ outside the support of
$\Psi'$. 

We will apply a result of Andreotti-Grauert \cite[p.109]{hk2} in order to solve
the equation $\dbar H=\Phi^*$ on $M'\times (N\setminus N'')$ with $H$
smooth of bidegree $(k-p,k-p-1)$. Let $\widetilde\chi$
be a cut-off function equal to 0 on $M'\times N''$ and 1 in a
neighbourhood of  $M'\times (N\setminus N')$. The form $\widetilde\chi H$ is
defined on $M'\times N$. It is clear that
$\Phi':=\Phi^*-\dbar(\widetilde\chi H)$ is horizontal in $M'\times
N'$ and satisfies $\dbar \Phi'=\Psi'$, which completes the proof.  

In order to apply the  Andreotti-Grauert theorem,
i.e. to solve the $\dbar$-equation for a $\dbar$-closed form of
bidegree $(l,k-s)$, $s\geq p$, in $M'\times (N\setminus \overline N'')$, we
only have to prove
that $M'\times (N\setminus \overline N'')$ satisfies the right convexity
property. More
precisely, one should construct
a smooth exhaustion function $\rho$ on  $M'\times (N\setminus
\overline N'')$ such
that $\ddc\rho$ has at every point $p+1$ strictly positive eigenvalues.
The domain is completely strictly $p$-convex
in the terminology of \cite[p.65]{hk2}. We need a much weaker result
than Theorem 12.7 in \cite{hk2}. 

Let $\rho_1$ be a smooth strictly convex function on $N$ such that
$\rho_1(z)\rightarrow\infty$ when $z\rightarrow \partial N$ and
$N''=\{ \rho_1 <1\}$. 
Since  $M'$ is strictly convex, we may find an unbounded exhaustion function $\rho_0$
for $M'$ which is smooth and strictly convex.
Define
\begin{equation*}
\rho(z)  := \rho_0(z')+c\rho_1(z'')+\kappa(\rho_1(z'')),\qquad
z=(z',z'')\in M'\times (N\setminus \overline N'').
\end{equation*}
with  $\kappa(t):=\frac{1}{t-1}$ and $c>0$ large enough.
The function $\rho$  is an exhaustion function   on   $M'\times
(N\setminus \overline N'')$.
The $p$ eigenvalues of $\ddc\rho$ with respect the the variable $z'$
are strictly positive.
On the other hand, since  
\begin{equation*}
i\partial \overline{\partial} (\kappa \circ \rho_1) =
\kappa'\cdot i \partial \overline{\partial}\rho_1+\kappa''\cdot i \partial \rho_1\wedge \overline{\partial} \rho_1,
\end{equation*}
and $\kappa''(t)\gg \vert \kappa'(t)\vert$ as $t\to 1^{+},$   
$\ddc\rho$ admits, at  every
point, at least one  strictly  positive eigenvalue with respect to the
variable $z''$. This completes the proof. 
\endproof 
 
\noindent
{\bf End of the proof of Theorem \ref{th_ddc}.} Define
$\Phi'':=-i\pi(\Phi'-\overline\Phi')$. 
This is a real smooth horizontal form in $M'\times N'$. 
We have
$$\ddc\Phi''=\ddbar (\Phi'-\overline\Phi')=\partial\Psi'+\dbar \overline\Psi'=d\Psi=\Omega'.$$
The smooth form $\Phi''$ is not necessarily positive. We can assume that it
has support in $M'\times F$ for some compact subset $F$ of $N'$. We now
construct a horizontal closed form $U$ on $M'\times N'$ of bidegree
$(k-p,k-p)$ which is strictly positive on $M'\times F$. Then, the form
$\Phi_2:=\Phi''+cU$, with $c>0$ large enough, is positive and
satisfies $\ddc \Phi_2=\Omega'$. 

For every point $z\in \overline M'\times F$ there is a complex plane $P$
of dimension $p$ passing through $z$ which does not intersect
$\overline M'\times\partial N'$. This plane defines by integration a
positive closed $(k-p,k-p)$-current $[P]$. Using a
convolution, we obtain by averaging on 
small pertubations of $[P]$, a smooth positive closed form $U_z$ which is
horizontal in $M'\times N'$ and is strictly positive at $z$. By
continuity, such a
form is strictly positive in a neighbourhood of $z$. It is enough to
take a finite sum of such forms in order to obtain a form $U$ which
is strictly positive on $\overline M'\times F$. This completes the proof.
\hfill $\square$

\begin{remark} \label{rk_ddbar_smooth}
\rm
If $\Omega$ is a continuous form then $\|\Phi\|_{\Cc^1(M'\times N')}\leq
c\|\Omega\|_{\Cc^0(M\times N'')}$ with a constant $c>0$ independent of
$\Omega$. Indeed, we are using a solution given by
a ``good'' kernel.
\end{remark}


\section{Horizontal-like maps} \label{section_horizontal}

In this section we introduce the class of horizontal-like maps,
 the main dynamical objects of our study, and we 
give some basic properties.

\medskip

\noindent
{\bf $\bullet$ Horizontal-like maps and Julia sets.} A horizontal-like map $f$ on $D$ 
is not necessarily defined on the whole domain $D$ but only on 
a vertical subset $f^{-1}(D)$ of $D$. It takes values in a horizontal subset 
$f(D)$ of $D$.  
Horizontal-like maps are defined by their graphs $\Gamma$ as follows 
\cite{ds6}. Let $\pr_1$ and $\pr_2$ be the canonical 
projections of $D\times D$ on its factors.

\begin{definition} \rm
 \label{def_horiz_map}
A {\it horizontal-like map} $f$ on $D$ is a holomorphic map with graph 
$\Gamma$ such that
\begin{enumerate}
\item $\Gamma$ is a submanifold of $D\times D$.
\item $\pr_{1|\Gamma}$ is injective; $\pr_{2|\Gamma}$ has finite fibers.
\item $\overline\Gamma$ does not intersect 
$\overline{\partial_v D} \times \overline D$ nor 
$\overline D \times \overline{\partial_h D}$.
\end{enumerate} 
\end{definition}
The last property is equivalent to the fact that the projections of
$\Gamma$ on the first factor $M$ and the last factor $N$ in $D\times
D$ are relatively compact.
The map $f=\pr_2\circ(\pr_{1|\Gamma})^{-1}$ is defined on $f^{-1}(D):=\pr_1(\Gamma)$ 
and its image is equal to $f(D):=\pr_2(\Gamma)$.
There exist open sets 
$M'\Subset M$ and $N'\Subset N$ such that
$f^{-1}(D)\subset D_v:=M'\times N$ and $f(D)\subset D_h:=M\times
N'$. We have $\Gamma\subset D_v\times D_h$. This property
characterizes horizontal-like maps and we often use it in order to
check that a map is horizontal-like. 
Since $\Gamma$ is a submanifold of $D\times D$, 
when $z$ tends to $\partial f^{-1}(D)\cap D$, 
$f(z)$ tends to $\partial_vD$.
When $z$ tends to $ \partial f(D)\cap D$, $f^{-1}(z)$ tends to 
$\partial_hD$.
So, the vertical part of $\partial f^{-1}(D)$ is sent into the
vertical part of $\partial f(D)$.
If $g$ is another horizontal-like map on $D$, $f\circ g$ is also a horizontal-like map.
When $p=k$, we obtain the polynomial-like maps studied in 
\cite{ds1}.

If $\pr_{2|\Gamma}$ is injective, we say that $f$ is {\it invertible}. 
In this case, up to a coordinate change (an exchange of horizontal and
vertical directions), 
$f^{-1}:\pr_2(\Gamma)\rightarrow\pr_1(\Gamma)$ is a horizontal-like map.
When $k=2$ and $p=1$,  we obtain the H{\'e}non-like maps \cite{duj, dds}.
In order to simplify the paper, we consider only invertible horizontal-like maps.

Small pertubations of an invertible horizontal map are still 
horizontal and invertible if one shrinks
slightly the domain $D$. Therefore, it is easy to construct large families of
such maps.

Define $f^n:=f\circ\cdots\circ f$ ($n$ times) {\it the iterate of order $n$
of $f$} and  $f^{-n}:=f^{-1}\circ\cdots\circ f^{-1}$ ($n$ times) its
inverse. Let $\Kc_+$ (resp. $\Kc_-$) denote the set of points $z\in D$
such that $f^n$ (resp. $f^{-n}$) are defined at $z$ for every $n\geq
0$. In other words, we have $\Kc_+:=\cap_{n\geq 0} f^{-n}(D)$ and  
$\Kc_-:=\cap_{n\geq 0} f^n(D)$. It is easy to check that $\Kc_\pm$ are
closed in $D$; $\Kc_+$ is
vertical and  $\Kc_-$ is
horizontal. We call $\Kc_+$ {\it the filled Julia set of $f$}
and  $\Kc_-$ {\it the filled Julia set of $f^{-1}$}. Their boundaries
are called {\it Julia sets}. Define also $\Kc:=\Kc_+\cap\Kc_-$. This
is a compact subset of $D$. We have $f^{-1}(\Kc_+)=\Kc_+$,
$f(\Kc_-)=\Kc_-$ and  $f^{\pm 1}(\Kc)=\Kc$, see \cite{ds6}.

\medskip

\noindent
{\bf $\bullet$ Dynamical degrees, Green currents and equilibrium measure.}
The operator $f_*:=(\pr_{2|\Gamma})_*\circ (\pr_{1|\Gamma})^*$ 
acts continuously on horizontal currents. 
If $S$ is a horizontal current or form, so is 
$f_*(S)$.
The operator $f^*:=(\pr_{1|\Gamma})_*\circ (\pr_{2|\Gamma})^*$ 
acts continuously on vertical currents. 
If $R$ is a vertical current or form, so is $f^*(R)$.
The  continuity of  $f^{\ast},$  $f_{\ast}$  for non-invertible maps is  treated in \cite{ds7}.  
Recall from \cite{ds6} the following proposition for positive closed
currents of the right bidegree.

\begin{proposition} \label{operators_on_currents} The operator $f_{\ast}:\Cc_h(D_v)\rightarrow \Cc_h(D_h)$
is well-defined and continuous. Moreover, there exists an integer $d\geq 1$ such that
$\| f_*(S)\|_h=d\| S\|_h$ for every $S\in \Cc_h(D_v)$. 
The operator $f^*:\Cc_v(D_h)\rightarrow \Cc_v(D_v)$
is well-defined and continuous. If $R$ belongs to $\Cc_v(D_h)$, we have 
$\| f^*(R)\|_v=d\| R\|_v$.
\end{proposition}
The integer $d$ is  called {\it the main  dynamical degree of
  $f$}. In the sequel, it is often denoted by $d(f).$ 
Note that the previous proposition implies that $d(f)=d(f^{-1})$ and $d(f^n)=d^n$. Consider
a vertical subvariety $L$ of dimension $k-p$ in $D$. The projection $\pi_2:L\rightarrow
N$ defines a (ramified) covering. If $m$ is the degree of this
covering, the current $[L]$ has slice mass $m$. We deduce from the
previous proposition that $f^{-1}(L)$ is a vertical subvariety of
degree $md$. For $m=1$, we obtain that $d$ is an integer. There is an
analogous picture when we push forward a horizontal subvariety. 
Note also that the projection of $\Gamma$ onto the product of the
first factor $N$ with the second factor $M$ defines a (ramified)
covering of degree $d$.
The following results were proved in \cite{ds6}.

\begin{theorem} \label{th_cv_mu}
Let $f$ be an invertible horizontal-like map on $D=M\times N$,  
$d$ its main dynamical degree and $\Kc_\pm$, $\Kc$
the filled Julia sets as above. Let $R$ and $S$ be smooth forms
in $\Cc^1_v(D)$ and  $\Cc^1_h(D)$ respectively.
Then $d^{-n} (f^n)^*(R)$ (resp. $d^{-n} (f^n)_*(S)$)
converge to a current $T_+$ in $\Cc^1_v(D)$ (resp. $T_-$ in  
$\Cc^1_h(D)$) which does not depend on $R$ (resp. $S$) and $d^{-2n}
(f^n)^*(R)\wedge (f^n)_*(S)$ converge to the probability measure
$\mu:=T_+\wedge T_-$. The current $T_+$ (resp. $T_-$) is supported on
the Julia set $\partial \Kc_+$ (resp. $\partial \Kc_-$) and is
invariant under $d^{-1} f^*$ (resp. under $d^{-1}f_*$).
The measure $\mu$ is invariant under $f^*$, $f_*$ and is supported
on $\partial \Kc_+\cap\partial\Kc_-$.
\end{theorem}

The current $T_+$ (resp. $T_-$) is {\it the Green current} associated
to $f$ (resp. $f^{-1}$). The measure $\mu$ is called {\it the
  equilibrium measure} of $f$.

\begin{theorem}
With the notation of the previous theorem, the topological
entropy of $f$ on $\Kc$ is equal to $\log d$ and $\mu$ is a measure of
maximal entropy $\log d$.
\end{theorem}

The notion of entropy will be recalled in Section \ref{section_measure}.
We now introduce {\it the other dynamical degrees} of $f$.
Recall that the open sets $M'\Subset M$ and $N'\Subset N$ are chosen
so that $f^{-1}(D)\subset M'\times N$ and $f(D)\subset M\times
N'$. So, the restriction of $f$ to $M'\times N'$ is also horizontal-like. 
For every $0\leq  s\leq p,$ let
\begin{equation*} 
d^{+}_s=d_s(f) := \limsup\limits_{n\to\infty}\Big\lbrace
\sup\limits_{S} \| (f^n)_{\ast} 
S\|_{M'\times N}  \Big\rbrace^{1/n},
\end{equation*}
the supremum being taken over all positive closed horizontal currents
$S$ of bidegree $(k-s,k-s)$ on $D'=M'\times N'$
such that  $\| S\|_{D'}=1.$   
For every $0\leq s\leq k-p,$ define 
\begin{equation*} 
d^{-}_s=d_s(f^{-1}):=  \limsup\limits_{n\to\infty}\Big\lbrace
\sup\limits_{R}  \| (f^n)^{\ast} R\|_{M\times N'} \Big \rbrace^{1/n},
\end{equation*}
the supremum being taken over all positive closed vertical  currents
$R$ of bidegree $(k-s,k-s)$ on $D'=M'\times N'$
such that $\| R\|_{D'}=1.$   
In the sequel we will write  for short
\begin{equation*} 
\delta_{+}:=d^{+}_{p-1}\qquad\text{and}\qquad \delta_{-}:=d^{-}_{k-p-1}.
\end{equation*}
These are the dynamical degrees which have to be compared to $d$.

\begin{lemma} \label{lemma_def_degree}
The dynamical degrees do not depend  on the choice of the particular convex domains
$M'$ and $N'$. Moreover, we have $d_0^+=d_0^-=1$ and $d^{+}_p=d^{-}_{k-p}=d$.
\end{lemma}
\proof
Let $M''$ and $N''$ be convex open sets such that $M''\Subset
M'\Subset M$, $N''\Subset N'\Subset N$ and $f^{-1}(D)\subset M''\times
N$, $f(D)\subset M\times N''$. If in the previous
definition, we replace $M'$ by $M''$  and $N'$ by $N''$, we obtain
$\delta_s^+$ and $\delta_s^-$. It is enough to prove that
$\delta_s^+= d_s^+$ and $\delta_s^-= d_s^-$. 
We prove the first equality; the second one is obtained in the same way.
Let $\widetilde S$ be a horizontal positive closed current of bidegree
$(k-s,k-s)$ on $M''\times N''$. Since $f$ is horizontal-like,
$f_*(\widetilde S)$ is horizontal in $M\times N''$ and there is a
constant $A>0$ independent of $\widetilde S$ such that 
$\|f_*(\widetilde S)\|_{M'\times N'}\leq A\|\widetilde S\|_{M''\times N''}$.
In particular, if $S$ is horizontal in $M'\times N'$
then we have $\|(f^n)_*S\|_{M'\times N'}\leq A\|(f^{n-1})_*S\|_{M''\times
  N''}$ for $n\geq 2$. If, moreover, $\|S\|_{M'\times N'}=1$, then $S':=f_*(S)$ is horizontal in
$M''\times N''$ with bounded mass. Therefore,
\begin{eqnarray*} 
d^{+}_s & = & \limsup\limits_{n\to\infty}\Big\lbrace
\sup\limits_{S} \| (f^n)_{\ast} 
S\|_{M'\times N'}  \Big\rbrace^{1/n} \\
& \leq &  \limsup\limits_{n\to\infty}\Big\lbrace
\sup\limits_{S} \| (f^{n-1})_{\ast}  S\|_{M''\times N''}
\Big\rbrace^{1/n} \\
& = &  \limsup\limits_{n\to\infty}\Big\lbrace
\sup\limits_{S} \| (f^{n-2})_{\ast}  S'\|_{M''\times N''}
\Big\rbrace^{1/n} \leq \delta_s^+.
\end{eqnarray*}        
For $S$ horizontal in $M''\times N''$ with
$\|S\|_{M''\times N''}=1$, define  also
$S':=f_*(S)$. Then $S'$ is horizontal in $M'\times N'$ with 
bounded mass and we have
\begin{eqnarray*} 
\delta^{+}_s & = & \limsup\limits_{n\to\infty}\Big\lbrace
\sup\limits_{S} \| (f^n)_{\ast} 
S\|_{M''\times N''}  \Big\rbrace^{1/n} \\
& \leq &  \limsup\limits_{n\to\infty}\Big\lbrace
\sup\limits_{S} \| (f^{n})_{\ast}  S\|_{M'\times N'}  \Big\rbrace^{1/n}\\
& = &  \limsup\limits_{n\to\infty}\Big\lbrace
\sup\limits_{S} \| (f^{n-1})_{\ast}  S'\|_{M'\times N'}
\Big\rbrace^{1/n} \leq d_s^+.
\end{eqnarray*}        
This implies the first part of the lemma. 

Since $f_*$ preserves the
mass of positive measures on $f^{-1}(D)$, we obtain that $d_0^+\leq
1$.  If $S$ is a probability measure on $\Kc$ then $(f^n)_*(S)$ is 
also a probability on $\Kc$. So, $d_0^+=1$. We obtain in the same way
that $d_0^-=1$.  

Assume that $S$ is of bidegree $(k-p,k-p)$. By definition of slices, 
we have $\|S\|_h\lesssim\|S\|_{M'\times N}$ and
as we already discussed in Section \ref{section_ddc}, 
$\|S\|_{M'\times N}\lesssim \|S\|_h$. So,
$$d^n\lesssim\|(f^n)_*S\|_{M'\times N}\lesssim d^n,$$
which implies that $d_p^+=d$. We obtain in the  same way that $d_{k-p}^-=d$. 
\endproof

\medskip

\noindent
{\bf $\bullet$ Action on super-functions.} We reduce slightly $D$ and
assume that $f$ is defined in a neighbourhood of $\overline D$.  Let
$\Phi$ be a current in $\PSH_h(D)$ and $\Lambda$ the super-function
associated to $\Phi$ defined on $\Cc_v^1(M\times \overline N)$, i.e. 
$\Lambda(R):=\langle R,\Phi\rangle$, 
see Proposition \ref{prop_super_function}. 
The following lemma is useful in our calculus.

\begin{lemma} \label{lemma_calculus}
The function $R\mapsto \Lambda\big(d^{-1}f^*(R)\big)$ is the
  super-function associated to $d^{-1} f_*(\Phi)$. In other words, we
  have 
$$\langle f^*(R),\Phi\rangle =\langle R,f_*(\Phi)\rangle$$ 
for  $R\in \Cc_v(D)$ and $\Phi\in \PSH_h(D)$. 
\end{lemma}
\proof
Let $\Lambda'$ denote the function  $R\mapsto
\Lambda\big(d^{-1}f^*(R)\big)$ and $\Lambda''$ the super-function
associated to  $d^{-1} f_*(\Phi)$. It is clear that
$\Lambda'(R)=\Lambda''(R)$ for $R$ smooth. 
We have to prove this equality for general $R$.

Let $\Rc$ be the current in
$\Delta\times D$ associated to the structural disc $(R_\theta)$
constructed in Section \ref{section_ddc}. If $F :\Delta\times
f^{-1}(D)\rightarrow\Delta\times f(D)$ is the map given by
$F(\theta,z):=(\theta,f(z))$, one can check that the current $d^{-1}F^*(\Rc)$ defines a
structural disc $(R_\theta')$ with
$R_\theta'=d^{-1}f^*(R_\theta)$. Since $\Lambda$ is p.s.h.,
$\Lambda'(R_\theta)=\Lambda(R'_\theta)$ is subharmonic on
$\theta\in\Delta$. The super-function $\Lambda''$ is also subharmonic
on the disc $(R_\theta)$ and coincide with $\Lambda'$ 
at $R_\theta$ with $\theta\not=1$ because $R_\theta$ is smooth for $\theta\not=1$. Hence, $\Lambda'$ and $\Lambda''$
coincide also at $R_1=R$, that is, $\Lambda'(R)=\Lambda''(R)$.
\endproof

\medskip

\noindent
{\bf $\bullet$ Product maps.} 
Let $f_i$ be horizontal-like maps on $D_i=M_i\times N_i$. Define  
the product map $F(x_1,x_2):=(f_1(x_1),f_2(x_2))$ on $D_1\times
D_2$. 
Up to a permutation of coordinates, we can identify $D_1\times D_2$
to $(M_1\times M_2)\times (N_1\times N_2)$. One checks easily that
$F$ is a
horizontal-like map on this domain. If $d_i$ denote the main dynamical degree of
$f_i$, the main degree of $F$ is $d_1d_2$. We can deduce from Theorem
\ref{th_cv_mu} the following properties.  
If $T_{i,\pm}$ are the Green currents associated to $f_i^{\pm 1}$,
the Green currents associated to $F^{\pm 1}$ are $T_{1,+}\otimes
T_{2,+}$ and $T_{1,-}\otimes T_{2,-}$. 
If $\mu_i$ are the equilibrium measures of $f_i$, the equilibrium
measure of $F$ is $\mu_1\otimes \mu_2$.
In what follows, we will
use the product $F$ of the horizontal-like maps $f_1:=f$ and
$f_2:=f^{-1}$ defined on $D=M\times N$ as above.
In this case, we have $M_1=N_2=M$ and $M_2=N_1=N$; the Green currents
of $F$ and $F^{-1}$ are $T_+\otimes T_-$ and $T_-\otimes T_+$.
We can perturb $F$ in order to obtain new families of examples.

\medskip

\noindent
{\bf $\bullet$ About the hypothesis on dynamical degrees.} 
The hypothesis we need in this paper is that the main
dynamical degree is larger than the other dynamical degrees.
The following proposition shows that the family of the maps $f$ satisfying
this condition is open.

\begin{proposition} \label{prop_pertub}
Let $f$ be a horizontal-like map on $D=M\times N$ with the main dynamical
degree $d$ as above
and 
$D':=M'\times N'$  a domain such that $D'\Subset
D$ and that $D\setminus D'$ is small enough. Then every small pertubation 
$f_\epsilon$ of $f$  is a
horizontal-like map on $D'$ of the same main dynamical degree $d$. 
If the dynamical degree of order $s$ of $f$ is strictly smaller than
$d$, then the dynamical degree of order $s$ of
$f_\epsilon$ satisfies the same property.
\end{proposition}
\proof
It is clear that $f_\epsilon$ is horizontal-like on $D'$. Since $d$
can be interpreted as the degree of a covering, the main
dynamical degree of $f_\epsilon$ is also $d$. Let $d^+_s$ and $d_{\epsilon,s}^+$ denote the
dynamical degrees of order $s$ of $f$ and $f_\epsilon$. 
Fix a constant $\delta$ such that $d_s^+<\delta<d$ and a domain
$D''=M''\times N''$ in $D'$ such that $D'\setminus D''$ is small
enough. So, $f$ and $f_\epsilon$ restricted to $D''$ are horizontal-like.
Consider a
horizontal positive closed $(k-s,k-s)$-current $S$ of mass 1 in $D''$. By
Lemma \ref{lemma_def_degree}, there is an integer $n_0$ independent of $S$ such that the mass of
$(f^{n_0})_*S$ on $D'$ is smaller than
$\delta^{n_0}/2$. If $\omega$ denotes the standard K{\"a}hler form on
$\C^k$, we have since $S$ is supported on $D''$ and
$f^{-n_0}(D'')\subset M''\times N$
$$\|(f^{n_0})_*S\|_{D''}=\int_{D''} (f^{n_0})_*S\wedge
\omega^s=\int_{f^{-n_0}(D'')\cap D''} S\wedge
(f^{n_0})^*\omega^s.$$
If $f_\epsilon$ is close enough to $f$, 
$(f_\epsilon^{n_0})^*\omega^s-(f^{n_0})^*\omega^s$ is
a small form on $f^{-n_0}(D')\cap D'$ and  $f^{-n_0}_\epsilon(D'')\cap D''\subset
f^{-n_0}(D')\cap D'$. Hence,
$$\|(f_\epsilon^{n_0})_*S\|_{D''}\leq \|(f^{n_0})_*S\|_{D'}+\int_{f_\epsilon^{-n_0}(D'')\cap D''} S\wedge
\big[(f_\epsilon^{n_0})^*\omega^s-(f^{n_0})^*\omega^s\big].$$
It follows that $\|(f_\epsilon^{n_0})_*S\|_{D''}\leq
\delta^{n_0}$. The estimate is independent of $S$ and implies by
iteration that $\|(f_\epsilon^n)_*S\|_{D''}\lesssim
\delta^n$ for $n\geq 1$ uniformly on $S$. Hence,
$d_{\epsilon,s}^+\leq \delta<d$. We get a similar results for
$f_\epsilon^{-1}$ and its dynamical degrees. 
\endproof


\section{Convergence theorems} \label{section_cv}

In this section we will give several quantitative versions of 
Theorem \ref{th_cv_mu} under the hypothesis that the main
dynamical degree $d$ is strictly larger than the degrees $\delta_+:=d^+_{p-1}$
and $\delta_-:=d^-_{k-p-1}$. We will see that this hypothesis is
natural and is satisfied for large families of maps. A similar
condition was considered in the context of polynomial-like maps, see
\cite{ds1}.

\medskip
\noindent
{\bf $\bullet$ Convergence towards the Green currents.}
We will use the PSH horizontal currents as test ``forms''. The above
solution of the $\ddc$-equation allows to write such a test current as
the sum of a PSH current with good estimates and a $\ddc$-closed
one. We obtain in particular the following result.

\begin{theorem} \label{th_pc}
Let $f$ be an invertible horizontal-like map on $D=M\times N$ and $d$,
$\delta_+$ its dynamical degrees as above. Assume that $d>\delta_+$.
Then the Green current $T_+$ of $f$ is PC.  
\end{theorem}

We first consider the $\ddc$-closed test currents.
The following result shows that in this case, 
without any hypothesis on
the dynamical degrees,
the convergence is exponentially fast and uniform. 

\begin{proposition} \label{prop_test_harmonic}
Let  $\Hc$ be a compact family of currents in $\Hc_h(D)$. Then there are
constants $A_0>0$ and $\lambda_0>1$ such that 
$$|\langle d^{-n}(f^n)^*R-T_+,\Psi\rangle|\leq A_0\lambda_0^{-n}$$
for all $R\in \Cc_v^1(M'\times N)$, $\Psi\in\Hc$ and $n\geq 0$. 
\end{proposition}
\proof
Reducing $D$ allows to assume that $R$ is in $\Cc_v^1(M'\times
\widetilde N)$ and $\Hc$ is compact in $\Hc_h(\widetilde M\times N')$.
There is a
constant $A'>0$ such that $|\langle d^{-n}(f^n)^*R,\Psi\rangle|\leq A'$
for all $R\in \Cc_v^1(D)$, $\Psi\in\Hc$ and $n\geq 0$. 
This follows from Proposition \ref{prop_test_closed} since
$(R,\Psi)\mapsto \langle
R,\Psi\rangle$ is continuous.
If $\Psi'$ is in $\Cc^1_h(\widetilde M\times N')$, we have and
$\langle T_+,\Psi'\rangle=1$ and $\langle d^{-n}(f^n)^*R,\Psi'\rangle=1$ for
every $R\in\Cc_v^1(D)$.
By adding to $\Psi$ a multiple of $\Psi'$,
we can assume that $\langle T_+,\Psi\rangle =0$ and we only need to
prove the estimate under this assumption. 
Assume also for simplicity that $A'=1$. 

Denote by $\Lambda_\Psi$ the super-function $\Lambda_\Psi(R):=\langle
R,\Psi\rangle$ 
and $L:=d^{-1}f^*$ the linear operator from $\Cc_v^1(D)$ into
$\Cc_v^1(M'\times\widetilde N)$.
Since $T_+$ is invariant, we have $\Lambda_\Psi\circ L^n(T_+)=0$.
Let $\Fc$ denote the set of pluriharmonic super-functions
  $\Lambda$ on $\Cc_v^1(M\times \overline N)$ such that $\Lambda(T_+)=0$ and
  $\|\Lambda\|_\infty \leq 1$.
Then, by Lemma \ref{lemma_calculus} and the assumption that $A'=1$,  $\Lambda_\Psi\circ L^n$ belongs
to $\Fc$ for $n\geq 1$ and we  have 
$$\langle d^{-n}(f^n)^*R-T_+,\Psi\rangle = \Lambda_\Psi\circ L^n (R).$$
So, by induction, it is enough to show that $\|\Lambda\circ L\|_\infty\leq
1/\lambda_0$ for $\Lambda$ in $\Fc$ and 
for some constant $\lambda_0>1$. 

Assume that no constant $\lambda_0$ satisfies the above
condition. Then there are $\Lambda\in \Fc$ and  
$R'\in \Cc_v^1(M'\times \widetilde N)$ such that $|\Lambda(R')|$ is
as close to 1 as we want. Recall that as in Section \ref{section_ddc} we can
construct a structural disc $\tau'$ (resp. $\tau$) such that
$\tau'(1)=R'$ (resp. $\tau(1)=T_+$). Moreover, 
$\tau'(0)$, $\tau(0)$ are equal to a fixed current $R_0$. These discs are parametrized by a
fixed neighbourhood $\Delta$ of $[0,1]$. By Harnack's inequality
applied to the non-vanishing harmonic function
$1-\Lambda\circ\tau'$ on $\Delta$, $|\Lambda(R_0)|$ is close to 1. Applying again
the Harnack's inequality to $1-\Lambda\circ \tau$, we deduce that
$|\Lambda(T_+)|$ is close to 1.
This contradicts the definition of $\Fc$.  
\endproof

\

\noindent
{\bf Proof of Theorem \ref{th_pc}.} Fix a constant $\delta$ such that 
$\delta_+<\delta<d$. Consider a test current $\Phi$ in a fixed compact
set of $\PSH_h(D)$.  Define $\Omega_0:=\ddc\Phi$ and
$\Omega_n:=(f^n)_*\Omega_0$. The currents $\Omega_n$ are positive of
bidegree $(k-p+1,k-p+1)$ and by definition of $\delta_+$, we have $\|\Omega_n\|_{M'\times N}\leq
A\delta^n\|\Omega_0\|$ with $A>0$ independent of $\Phi$.
By Theorem \ref{th_ddc} applied to $M''$ and $M'$, there are negative horizontal $L^1$ forms $\Phi_n$ such that
$\ddc\Phi_n=\Omega_n$ with $\|\Phi_n\|_{M''\times N}\lesssim
\delta^n$. Then, $\delta^{-n}\Phi_n$ belong to a fixed compact set of
$\PSH_h(M''\times N)$. Define $\Psi_0:=\Phi-\Phi_0$ and
$\Psi_n:=f_*(\Phi_{n-1})-\Phi_n$ for $n\geq 1$.
We have $\ddc \Psi_n=0$ and since $f_*$ is continuous,
$\delta^{-n}\Psi_n$ belong to some compact set in 
 $\Hc_h(M''\times N)$.

Fix a current $R$ in $\Cc_v^1(D)$. We can assume that $M''$ is chosen
so that $R$ is supported on $M''\times N$.
We have since $\Phi=\Psi_0+\Phi_0$
\begin{eqnarray*}
\langle d^{-n}(f^n)^*R,\Phi\rangle &=& \langle d^{-n}
(f^n)^*R,\Psi_0\rangle + \langle d^{-n+1} (f^{n-1})^*R, d^{-1}
f_*(\Phi_0)\rangle \\
&=& \langle d^{-n}
(f^n)^*R,\Psi_0\rangle + \langle d^{-n+1} (f^{n-1})^*R, d^{-1}
\Psi_1\rangle  \\
& & + \langle d^{-n+1} (f^{n-1})^*R, d^{-1} \Phi_1\rangle.
\end{eqnarray*}
By induction and using the identity $f_*(\Phi_n)=\Psi_{n+1}+\Phi_{n+1}$, we obtain
\begin{eqnarray} \label{eq_current_iterate}
\langle d^{-n}(f^n)^*R,\Phi\rangle &=&  \sum_{0\leq j\leq n-1}
\langle d^{-n+j} (f^{n-j})^*R,d^{-j}\Psi_j\rangle +  \langle R,
d^{-n}f_*(\Phi_{n-1})\rangle \nonumber\\
& = & \sum_{0\leq j\leq n}
\langle d^{-n+j} (f^{n-j})^*R,d^{-j}\Psi_j\rangle + 
 \langle R, d^{-n}\Phi_n\rangle.
\end{eqnarray}
Now assume that $R$ is smooth and let $n\rightarrow\infty$. 
The estimate on $\|\Phi_n\|$ implies that the last term tends to
0. Recall that $\delta^{-n}\Psi_n$ belong to a compact set in
$\Hc_h(M''\times N)$ and that $\delta <d$. 
On the other hand,
by Theorems \ref{th_cv_mu}, $d^{-n+j} (f^{n-j})^*R$ tends
to $T_+$ when $n-j\rightarrow\infty$.
Proposition \ref{prop_test_harmonic} and Lebesgue's convergence theorem, applied to the
series in the identity (\ref{eq_current_iterate}), imply that for
$\Phi$ smooth
 \begin{eqnarray} \label{eq_current_limit}
\langle T_+,\Phi\rangle &=&  \Big\langle T_+,\sum_{j\geq 0} d^{-j} \Psi_j \Big\rangle.
\end{eqnarray}
Observe that the last sum is pluriharmonic and depends continuously on
$\Phi$ in $\PSH_h(D)$. 
It follows from Proposition \ref{prop_test_closed} that the right hand side of the
last identity depends continuously on $\Phi$. So, $T_+$ is a PC
current and the identity (\ref{eq_current_limit}) holds for all $\Phi$
in $\PSH_h(D)$.
\hfill $\square$
\\

The following propositions give the speed of convergence towards the
Green current.

\begin{proposition} \label{prop_current_speed_above}
Let $f$ be as in Theorem \ref{th_pc} with $d>\delta_+$.
Let $\Pc_v$ be a compact family of currents in $\Cc_v^1(D)$ and $\Dc_h$ a
compact family of test currents in $\PSH_h(D)$.  
Then, there exist constants $A>0$ and
 $\lambda>1$ such that 
$$\langle d^{-n} (f^n)^*R-T_+,\Phi\rangle\leq A\lambda^{-n}$$
for all $R\in\Pc_v$, $\Phi\in \Dc_h$ and $n\geq 0$.
\end{proposition}
\proof
Observe that when $\Phi$ belongs to a compact family  in $\PSH_h(D)$,
$\delta^{-n}\Phi_n$ and $\delta^{-n}\Psi_n$  belong to compact families in $\PSH_h(M'\times N)$ and 
in $\Hc_h(M'\times N)$ for some $M'\Subset M$. 
It follows from the identities (\ref{eq_current_iterate}) and
(\ref{eq_current_limit}) that $\langle d^{-n}
(f^n)^*R-T_+,\Phi\rangle$ is equal to 
$$\sum_{0\leq j\leq n} d^{-j}
\big\langle d^{-n+j} (f^{n-j})^*R-T_+,\Psi_j \big\rangle  - 
\sum_{j\geq n+1} d^{-j} \langle T_+, \Psi_j \rangle +  d^{-n}\langle
R, \Phi_n\rangle.$$
Proposition \ref{prop_test_harmonic} implies that $|\langle d^{-n+j}
(f^{n-j})^*R-T_+,\Psi_j\rangle|
\lesssim \lambda_0^{-n+j}\delta^j$.
We also deduce from Proposition \ref{prop_test_closed} applied to
$M'\times N$ instead of $D$, that $|\langle
T_+,\Psi_j\rangle|\lesssim \delta^j$.
Since $\Phi_n$ is negative, the last term in the previous sum is negative. 
This implies the desired estimate  for $1<\lambda<\min(\lambda_0, d/\delta)$.
\endproof

\begin{proposition} \label{prop_current_speed}
Let $f$ be as in Theorem \ref{th_pc} with $d>\delta_+$.
Let $\Pc'_v$ be a bounded family of PB currents in $\Cc_v^1(D)$ 
and $\Dc_h$ a compact family of test currents in $\PSH_h(D)$. Then,
there exist constants $A>0$ and 
 $\lambda>1$ such that  
$$|\langle d^{-n} (f^n)^*R-T_+,\Phi\rangle|\leq A\lambda^{-n}$$
for all $R\in\Pc'_v$, $\Phi\in \Dc_h$ and $n\geq 0$.
\end{proposition}
\proof
 As in Proposition \ref{prop_current_speed_above}, it is enough
to estimate $|\langle
R, \Phi_n\rangle|$.
We have
$|\langle R, \Phi_n\rangle|\lesssim \delta^n$
since $R$ belongs to a bounded family of PB currents in
$\Cc^1_v(M'\times N)$ for some $M'\Subset M$.
This implies the proposition.
\endproof

\begin{remark} \label{rk_degree}
\rm
In Propositions \ref{prop_current_speed_above} and \ref{prop_current_speed}, the condition $d>\delta_+$ is
superflous if the mass of $d^{-n} (f^n)_*(\ddc\Phi)$ decreases to 0
exponentially and uniformly on $\Phi\in\Dc_h$ when $n$ goes to
infinity.
We will use this observation in the proof of Theorem \ref{th_mixing}.
\end{remark}

The following result gives a strong ergodic property for the action of $f$
on vertical currents.

\begin{theorem}
Let $f$ be an invertible horizontal-like map as above with $d>\delta_+$. Then
$d^{-n} (f^n)^*R$ converge to $T_+$ uniformly on $R\in\Cc_v^1(D)$. In
particular, $T_+$ is the unique current in $\Cc_v^1(D)$ which is
invariant under $d^{-1}f^*$.
\end{theorem}
\proof
Since smooth horizontal test forms are generated by the PSH ones,
it is enough to test smooth PSH horizontal forms. 
Using identity (\ref{eq_current_iterate}) for $\Phi$ smooth, we
only have to show
that $d^{-n}\langle R,\Phi_n\rangle$ tend to 0 uniformly on
$R$. Recall that $\Phi_n$ is negative, so $d^{-n}\langle
R,\Phi_n\rangle$ is negative. For
simplicity, we reduce the size of $D$ and we replace $R$ by $d^{-1}
f^*(R)$. So, we can assume that $f$ is
defined in a neighbourhood $\widehat D=\widehat M\times \widehat N$ 
of $\widetilde D=\widetilde M\times\widetilde N$ and that $R$, $\Phi_n$
are vertical or horizontal on $M''\times\widehat N$ and $\widehat
M\times N''$ respectively. Recall that the
convex sets 
$\widetilde M$, $\widehat M$, $\widetilde N$ and $\widehat N$ are chosen so that $M\Subset
\widetilde M\Subset \widehat M$ and $N\Subset \widetilde N\Subset
\widehat N$. We can also assume that the $\Cc^1$ norm of $f^{-1}$ on
$\widehat D$ is bounded by a constant $A>0$.

Assume by contradiction that there is an increasing
sequence $(n_i)$ such that $\langle R_i,\Phi_{n_i}\rangle\leq -2cd^{n_i}$ for
some positive constant $c>0$ and some sequence $(R_i)$ in
$\Cc_v^1(M'\times \widehat N)$. 
Let $(R_{i,\theta})_{\theta\in\Delta}$ denote the structural
discs associated to $R_i$ as in Section \ref{section_ddc}. 
Define $\varphi_i(\theta):=\delta^{-n_i}\langle
R_{i,\theta},\Phi_{n_i}\rangle$ with $\delta_+<\delta<d$. The properties of $R_{i,\theta}$ and of $\Phi_n$
imply that $\varphi_i$ belong to a compact family of subharmonic
functions on $\Delta$. It is then classical that for every compact subset $K$ of
$\Delta$ there are  constants $C>0$ and $\alpha>0$ such that
$\|e^{-\alpha\varphi_i}\|_{L^1(K)}\leq C$, see e.g. \cite{Hormander}. 

The currents 
$R_{i,\theta}$ are obtained by smoothing of $R$. Using a coordinate change, we obtain that
$$\langle R_{i,\theta},\Phi_{n_i}\rangle = \langle R_i,\Phi_{n_i,\theta}\rangle$$
where $\Phi_{n_i,\theta}$ is a smoothing of $\Phi_{n_i}$. With the
notation in Section \ref{section_ddc}, we have
$$\Phi_{n_i,\theta}:=\int (h_{a,b,\theta})^*(\Phi_{n_i}) d\nu(a,b).$$
 Since the family $h_{a,b,\theta}$ is
holomorphic and $h_{a,b,1}=\id$, we obtain
(see also \cite[Lemma 2.7]{ds6})
$$\|\Phi_{n_i,\theta}-\Phi_{n_i}\|_{\infty,D}\lesssim
|\theta-1|\|\Phi_{n_i}\|_{\Cc^1(\widetilde D)}$$
for $\theta$ close to 1. On the other hand,  the $\Cc^1$
norm of $f^{-n}$ is bounded by $A^n$, hence Theorem
\ref{th_ddc} and Remark \ref{rk_ddbar_smooth} imply that
$$\|\Phi_n\|_{\Cc^1(\widetilde D)}\lesssim
\|(f^n)_*(\ddc\Phi)\|_{\Cc^1(\widehat D)}\lesssim A^{2kn}.$$
Therefore, $\|\Phi_{n_i,\theta}-\Phi_{n_i}\|_{\infty,D}\lesssim
|\theta-1|A^{2kn_i}$ and since the mass of $R_i$ is bounded 
$$|\langle R_i, \Phi_{n_i,\theta}-\Phi_{n_i}\rangle|\lesssim |\theta-1|A^{2kn_i}.$$
Hence,  for $\theta$ in a disc of center 1 and of
radius $\simeq A^{-2kn_i}$, we have $\langle
R_{i,\theta},\Phi_{n_i}\rangle\leq -cd^{n_i}$ and then $\varphi_i(\theta) \leq
-cd^{n_i}\delta^{-n_i}$. This
contradicts the above uniform integrability of $e^{-\alpha\varphi_i}$.
\endproof

\medskip

\noindent
{\bf $\bullet$ Convergence towards the equilibrium measure.}
The main result in this section is the following property of the
equilibrium measure.

\begin{theorem} \label{th_pc_mu}
Let $f$ be an invertible horizontal-like map as above with
  $d>\delta_+$ and $d>\delta_-$. Then the equilibrium measure $\mu$ of
  $f$ is PC.
\end{theorem}
\proof
By Theorem \ref{th_pc}, $T_+$, $T_-$ are PC on $D$ and also on
$D':=M'\times N'$.
If $\varphi$ is a p.s.h. function on $D$,
$\varphi$ is locally integrable with respect to the trace measure
$T_-\wedge\omega^p$ of
$T_-$. Hence, $\varphi T_-$ defines a PSH horizontal
current. Moreover, the fact that $T_-$ is PC implies that $\varphi\mapsto\varphi T_-$
is continuous on $\varphi\in\PSH(D)$ with values in $\PSH_h(D)$.
Indeed, if $\Theta$ is a smooth vertical $(p,p)$-form, then
$\varphi\mapsto \langle \Theta,\varphi T_-\rangle$ is continuous,
since it is upper semi-continuous when $\Theta$ is positive and is
continuous when $\Theta$ is positive closed.
Using the PC property of $T_+$ and
the identity
$\langle\mu,\varphi\rangle=\langle T_+,\varphi T_-\rangle$, see
Remark \ref{rk_int_psh}, 
we obtain that $\langle \mu,\varphi\rangle$ depends continuously on
$\varphi$. Therefore, $\mu$ is PC.
\endproof

We can now prove estimates on the speed of convergence towards the
equilibrium measure.

\begin{proposition} \label{prop_measure_speed_above}
Let $f$ be as in Theorem \ref{th_pc_mu} with $d>\delta_+$ and $d>\delta_-$.
Let $\Pc_v$ (resp. $\Pc_h$) be a compact family of currents in
$\Cc_v^1(D)$ (resp. in $\Cc_h^1(D)$).
Then, there exist constants $A>0$ and
 $\lambda>1$ such that 
$$\langle d^{-2n} (f^n)^*R\wedge (f^n)_*S-\mu,\varphi\rangle\leq A\lambda^{-n}$$
for all $R\in\Pc_v$, $S\in\Pc_h$, $\varphi$ p.s.h. on $D$ with
$|\varphi|\leq 1$ and $n\geq 0$.
\end{proposition}
\proof
Since $\mu=T_+\wedge T_-$, we can write $\langle d^{-2n} (f^n)^*R\wedge (f^n)_*S-\mu,\varphi\rangle$
as the sum of the following two integrals
$$\langle d^{-2n} (f^n)^*R\wedge (f^n)_*S- d^{-n} (f^n)^*R\wedge
T_-,\varphi\rangle
=  \langle d^{-n} (f^n)_*S- T_-,\varphi d^{-n} (f^n)^*R\rangle$$
and
$$\langle d^{-n} (f^n)^*R\wedge T_-- T_+\wedge
T_-,\varphi\rangle
=  \langle d^{-n} (f^n)^*R- T_+,\varphi T_-\rangle.$$
Since $R$ is in a compact family in $\Cc_v^1(D)$, $d^{-n}
(f^n)^*R$ belong also  to a compact family in $\Cc_v^1(D)$ independent
of $n\geq 0$. Indeed, their supports are controlled.
Hence, for $|\varphi|\leq 1$, $\varphi d^{-n} (f^n)^*R$
belong to a compact family in $\PSH_v(D)$. By Proposition
\ref{prop_current_speed_above} applied to $f^{-1}$, the first integral is $\lesssim
\lambda^{-n}$ for some $\lambda>1$. Since $\varphi T_-$ belongs
to a compact family in $\PSH_h(D)$, the second integral is also
$\lesssim \lambda^{-n}$ for some $\lambda>1$. The proposition follows. 
\endproof

\begin{proposition} \label{prop_measure_speed}
Let $f$ be as in Theorem \ref{th_pc_mu} with $d>\delta_+$ and $d>\delta_-$.
Let $\Pc_v$ (resp. $\Pc_h$) be a bounded family of PB currents in
$\Cc_v^1(D)$ (resp. in $\Cc_h^1(D)$).
Then, there exist constants $A>0$ and
 $\lambda>1$ such that 
$$|\langle d^{-2n} (f^n)^*R\wedge (f^n)_*S-\mu,\varphi\rangle|\leq A\lambda^{-n}$$
for all $R\in\Pc_v$, $S\in\Pc_h$, $\varphi$ p.s.h. on $D$ with
$|\varphi|\leq 1$ and $n\geq 0$.
\end{proposition}
\proof
We proceed as in the proof of Proposition
\ref{prop_measure_speed_above} using 
Proposition \ref{prop_current_speed} 
instead of Proposition \ref{prop_current_speed_above}.
\endproof


\section{Properties of  the equilibrium measure} \label{section_measure}

In this section, we prove two important properties of the equilibrium
measure for horizontal-like maps with large main dynamical degree.

\medskip
\noindent
$\bullet$ {\bf Decay of correlations.} It was proved in \cite{ds6} that
the equilibrium measure is mixing for a general invertible
horizontal-like map. Under 
our hypothesis on dynamical degrees, we have the following
result.

\begin{theorem} \label{th_mixing}
Let $f$ be an invertible horizontal-like map as above with $d>\delta_+$ and
$d>\delta_-$. Then the equilibrium measure $\mu$ of $f$ is
exponentially mixing. More precisely, for all test functions $\phi$ of
class $\Cc^\alpha$ and $\psi$ of class $\Cc^\beta$ on $D$ with
$0<\alpha,\beta\leq 2$, the following estimate holds
$$|\langle \mu, (\phi\circ f^n)\psi\rangle -
\langle\mu,\phi\rangle\langle\mu,\psi\rangle |\leq A_{\alpha,\beta}
\lambda^{-n\alpha\beta}\|\phi\|_{\Cc^\alpha}\|\psi\|_{\Cc^\beta}$$
where $A_{\alpha,\beta}>0$ is a constant independent
of $\phi$, $\psi$, $n$ and $\lambda>1$ is a constant independent of
$\alpha$, $\beta$, $\phi$, $\psi$, $n$.   
\end{theorem}

Recall that the measure $\mu$ is mixing means that the left hand side of the above
inequality tends to 0 when $n$ goes to infinity. It follows from the theory of
interpolation between Banach spaces \cite{tri} that the previous inequality for
general $\alpha$, $\beta$ is deduced from the case where
$\alpha=\beta=2$, see \cite{din} for details. In the case of
H{\'e}non-like maps, i.e. $k=2$, we have $\delta_+=\delta_-=1$. So, the
hypothesis in the previous theorem is automatically
satisfied and we obtain the
following corollary.

\begin{corollary}
Let $f$ be a H{\'e}non-like map. Then the equilibrium measure of $f$ is
exponentially mixing.
\end{corollary}

\noindent
{\bf Proof of Theorem \ref{th_mixing}.} We only have to consider the
case where $\alpha=\beta=2$. 
Define
$$I_n(\phi,\psi):=\langle \mu, (\phi\circ f^n)\psi\rangle -
\langle\mu,\phi\rangle\langle\mu,\psi\rangle.$$
Observe that since $I_{n+1}(\phi,\psi)=I_n(\phi\circ f,\psi)$,
it is enough to consider the case where $n$ is even. 
Note also that since $\mu$ is invariant, $I_n(\phi,\psi)=0$ when
$\phi$ or $\psi$ is constant.

Near $\supp(\mu)$ we can write $\phi$ and $\psi$
as differences of functions which are strictly p.s.h. on a
neighbourhood of $\overline D$. So, we
can assume that $\ddc\phi\geq \ddc\|z\|^2$, $\ddc\psi\geq \ddc\|z\|^2$
and that $\phi$, $\psi$ have $\Cc^2$ norms bounded by a fixed
constant. This allows to fix a constant $A>0$ large enough such that
$(\phi(z)+A)(\psi(z')+A)$ and $(-\phi(z)+A)(\psi(z')-A)$ are
p.s.h. on $(z,z')$ in $D^2$. We have to bound from above 
$$I_{2n}(\phi,\psi)=I_{2n}(\phi+A,\psi+A)$$ 
and 
$$-I_{2n}(\phi,\psi)=I_{2n}(-\phi+A,\psi-A).$$
We will consider the first quantity, the proof for the second one is similar.
For that purpose, we will apply Proposition
\ref{prop_current_speed_above} and Remark \ref{rk_degree} to the
product $F$ of the horizontal-like maps $f$ and $f^{-1}$ defined in Section \ref{section_horizontal}.

Define $\varphi(z,z'):=(\phi(z)+A)(\psi(z')+A)$. Let $\Delta$ denote
the diagonal of $D\times D$ and $[\Delta]$ the current of integration
on $\Delta$. We have since $\mu$ is invariant
$$I_{2n}(\phi+A,\psi+A)  =  \big\langle \mu, (\phi\circ
f^n+A)(\psi\circ f^{-n}+A) \big\rangle - \langle
\mu,\phi+A\rangle\langle\mu,\psi+A\rangle.$$
Lifting these integrals to $D\times D$ and using the identity
$$d^{-2}F^*(T_+\otimes T_-)=d^2F_*(T_+\otimes
T_-)=T_+\otimes T_-,$$ 
we obtain that $I_{2n}(\phi+A,\psi+A)$ is equal to
\begin{eqnarray*}
\lefteqn{\langle (T_+\otimes T_-)\wedge [\Delta], \varphi\circ
F^n \rangle - \langle \mu\otimes\mu,\varphi\rangle} \\
 & = &  \big\langle (T_+\otimes T_-)\wedge  d^{-2n}(F^n)_*[\Delta], \varphi
  \big\rangle -  \big\langle (T_+\wedge T_-)\otimes (T_+\wedge T_-),\varphi \big\rangle \\
 & = &  \big\langle (T_+\otimes T_-)\wedge  d^{-2n}(F^n)_*[\Delta], \varphi
  \big\rangle -  \big\langle (T_+\otimes T_-)\wedge (T_-\otimes
 T_+),\varphi\big\rangle  \\
& = &   \big\langle d^{-2n}(F^n)_*[\Delta] - T_-\otimes T_+, \varphi
(T_+\otimes T_-)  \big\rangle. 
\end{eqnarray*} 
The current $[\Delta]$ is not horizontal but $F_*[\Delta]$ is
horizontal. So, we can apply Proposition
\ref{prop_current_speed_above} and Remark \ref{rk_degree} for $F^{-1}$.

For Remark \ref{rk_degree}, we need to show that 
the mass of 
$$d^{-2n} (F^n)^* \big[\ddc \varphi \wedge
(T_+\otimes T_-) \big]= \ddc(\varphi\circ F^n)\wedge (T_+\otimes T_-)$$ 
decreases exponentially (we reduce the size of $D$
if necessary). We have
$$\|\ddc(\varphi\circ F^n) \wedge (T_+\otimes T_-)\|_{D^2} = \int_{D^2}
 \big(\ddc\|z\|^2+\ddc\|z'\|^2 \big)^{k-1} \wedge \ddc (\varphi\circ F^n)\wedge (T_+\otimes T_-).$$
In the last wedge-product, $T_+$ depends only on $z$ and $T_-$ depends
only on $z'$. Then we expand 
$$\ddc(\varphi\circ
F^n)=\ddc \big[(\phi(f^n(z))+A)(\psi(f^{-n}(z'))+A) \big].$$ 
In this product, the terms containing 
mixed derivatives 
$$d\phi(f^n(z))\wedge \dc \psi(f^{-n}(z'))\quad \mbox {and} \quad
\dc\phi(f^n(z))\wedge d \psi(f^{-n}(z'))$$ 
vanish when wedged with 
$ \big(\ddc\|z\|^2+\ddc\|z'\|^2 \big)^{k-1} \wedge (T_+\otimes T_-)$
by bidegree consideration. 
This, combined with the fact that 
$\varphi$ and $\psi$ are bounded, implies that 
$$\|\ddc(\varphi\circ F^n) \wedge (T_+\otimes T_-)\|_{D^2}\lesssim \|\ddc(\phi\circ f^n)
\wedge T_+\|_D +\|\ddc(\psi\circ f^{-n}) \wedge T_-\|_D.$$
We have
$$\ddc(\phi\circ f^n) \wedge T_+ = d^{-n} (f^n)^*(\ddc\phi\wedge
T_+)$$ and
$$\ddc(\psi\circ f^{-n}) \wedge T_- = d^{-n} (f^n)_*(\ddc\psi\wedge
T_-).$$
Since $d>\delta_+$ and $d>\delta_-$, the masses of these currents 
decrease exponentially. This completes the proof.
\hfill $\square$ 

\begin{remark}\rm
We can prove the converse of Theorem \ref{th_pc}: {\it if the
current $T_+$ of $f$ is PB then $d>\delta_+$}. This will allow to prove
that $F$ satisfies also the hypothesis on dynamical degrees if
$d>\delta_+$ and $d>\delta_-$, hence we can apply directly
Proposition \ref{prop_current_speed_above}. However, the proof requires a
long development on the notion of 
super-functions introduced in Section \ref{section_ddc}, and we
prefer to avoid it here, see also \cite{ds8}. 
\end{remark}

\medskip
\noindent
$\bullet$ {\bf Lyapounov exponents.} We will show that when the main
dynamical degree of $f$ is larger than the other ones, the measure
$\mu$ is hyperbolic: it admits $p$ strictly positive and $k-p$
strictly negative Lyapounov exponents. We follow the approach by de
Th{\'e}lin \cite{deThelin}.

Recall that the measure $\mu$ is mixing and is supported on the
filled Julia set $\Kc:=\Kc_+\cap\Kc_-$ which is compact in $D$, see
\cite{ds6}. Using the theory of
Oseledec-Pesin \cite{pe}, we can decompose the tangent space of $\C^k$
at $\mu$-almost every point $x$ into a direct sum of vector subspaces
$T_x=\oplus_{i=1}^m E_{i,x}$ with the
following properties:
\begin{enumerate}
\item[-] The integer $m$ and the dimension of each $E_{i,x}$ do not
      depend on $x$.
\item[-] The decomposition $T_x=\oplus_{i=1}^m E_{i,x}$ is unique and depends in a
      measurable way on $x$.
\item[-] The vector bundle $E_{i,x}$ is {\it invariant} under $f$, that is, the
      differential $Df$ of $f$ defines an isomorphism between
      $E_{i,x}$ and $E_{i,f(x)}$.
\item[-] The decomposition $T_x=\oplus_{i=1}^m E_{i,x}$ has a
      {\it tempered distortion}. More precisely, if $I$ and $J$ are
      disjoint subsets of $\{1,\ldots,m\}$, define
      $E_{I,x}:=\oplus_{i\in I} E_{i,x}$ and  $E_{J,x}:=\oplus_{i\in
        J} E_{i,x}$. Then, 
the angle $\measuredangle \big(
      E_{I,f^n(x)}, E_{J,f^n(x)}\big)$ between
      $E_{I,f^n(x)}$ and $E_{J,f^n(x)}$ satisfies 
$$\lim_{n\rightarrow\pm\infty} {1\over n}
      \log \sin\measuredangle \big( E_{I,f^n(x)}, E_{J,f^n(x)}\big)=0.$$
\item[-] There are distinct real numbers $\lambda_i$ independent of $x$ such
      that 
$$\lim_{n\rightarrow\pm\infty}{1\over
        |n|} \log{\|Df^n(v)\|\over \|v\|}=\pm \lambda_i$$
uniformly on $v$ in $E_{i,x}\setminus\{0\}$.
\end{enumerate}
The constants $\lambda_i$ are {\it the Lyapounov exponents} of
$\mu$. {\it The multiplicity} of $\lambda_i$ is the dimension of
$E_{i,x}$. So, $\mu$ admits $k$ Lyapounov exponents counted with multiplicities\footnote{if $f$ is
  considered as a real map, the multiplicity of $\lambda_i$ is $2\dim
E_{i,x}$ and $\mu$ has $2k$ Lyapounov exponents; this is the reason for
the coefficients ${1\over 2}$ in Theorem \ref{th_hyperbolic}.}. 
The Lyapounov exponents of $f^n$ are $n\lambda_i$ even for $n$
negative. When there is no zero Lyapounov exponent, $\mu$ is said to
be {\it hyperbolic}.

\begin{theorem} \label{th_hyperbolic}
Let $f$ be an invertible horizontal-like map as above with dynamical degrees $d$,
$d_s^+$ and $d_s^-$. Define $\widetilde\delta_+:=\max_{s\leq p-1}
d_s^+$ and $\widetilde\delta_-:=\max_{s\leq k-p-1}d_s^-$. If
$\widetilde\delta_+<d$, then 
$\mu$ admits $p$ 
strictly positive Lyapounov exponents larger than or equal to ${1\over
  2k}\log(d/\widetilde\delta_+)$. If $\widetilde\delta_-<d$, then $\mu$
admits $k-p$ strictly negative ones which are smaller than or equal to $-{1\over
  2k}\log(d/\widetilde\delta_-)$.
\end{theorem}

We prove the first assertion. The second one is treated in the same
way using $f^{-1}$ instead of $f$. 
We will need the following lemmas where $\omega_v$ denotes the
restriction to $M''\times N$ of the standard
K{\"a}hler form $\omega$ on $\C^k$. Define $\widetilde d_q^+:=\max_{s\leq
  q} d_s^+$ for $0\leq q\leq p-1$. 

\begin{lemma}\label{lemma_pre_lov}
Let $\delta$ be a constant strictly larger than $\widetilde d_q^+$.
Then  there exists a constant $C>0$  such that for all positive closed
current $S$ of bidegree
$(k-q,k-q)$ supported on $M\times N'$
we have
\begin{equation*}
\int S\wedge (f^{n_1})^*\omega_v\wedge\ldots\wedge
(f^{n_q})^*\omega_v  
\leq C\delta^{n_1} \|S\|_D
\end{equation*}
for all integers $n_1\geq\cdots\geq n_q\geq 0$.
\end{lemma}
\proof We prove the lemma  by induction on $q$. Clearly, the lemma is  valid for $q=0$.
Suppose it holds for the rank $q-1$. This, applied to
$f$ restricted to $D_v:=M'\times N$ and to $S':=  (f^{n_q})_* S\wedge\omega$, implies that
\begin{equation*}
\int_{D_v} S'\wedge (f^{n_1-n_q})^*\omega_v\wedge\ldots\wedge (f^{n_{q-1} -n_q})^*\omega_v 
 \leq C\delta^{n_1-n_q} \|S'\|_{D_v}.
\end{equation*}
By definition of $d_s^+$, there is a constant $c>0$ such that
\begin{equation*}
\|S'\|_{D_v} = \| (f^{n_q})_* S \|_{D_v}\leq  c\delta^{n_q}\|S\|_D.
\end{equation*}
Consequently,
\begin{equation*}
\int_{D_v} S'\wedge (f^{n_1-n_q})^*\omega_v\wedge\ldots\wedge (f^{n_{q-1} -n_q})^*\omega _v
 \leq C\delta^{n_1} \| S\|_D.
\end{equation*}
for some constant $C>0$. 
The left hand side of the last inequality is equal to
$$\int_{f^{-n_q}(D_v)} S\wedge (f^{n_1})^*\omega_v\wedge\ldots\wedge
(f^{n_q})^*\omega_v.$$  
This implies  the lemma for rank $q$. Note that the last integral does
not change if we replace $f^{-n_q}(D_v)$ by $D$ since
  $(f^{n_q})^*\omega_v$ is supported on $f^{-n_q}(D_v)$.
\endproof

Let $\Gamma_n$ denote the graph of $(f,\ldots, f^n)$, i.e.
the set of points $(x,f(x),\ldots, f^n(x))$ in $D^{n+1}$. We will use
the standard K{\"a}hler metric $\omega_n$ in $D^{n+1}\subset \C^{k(n+1)}$.
If $\Pi_j$, with $0\leq j\leq n$, denote the projections from $D^{n+1}$ onto its
factors $D$, we have $\omega_n=\sum \Pi_j^*(\omega)$.
Let $\pi_j$ denote the restriction of $\Pi_j$ to 
$\Gamma_n$ and $\volume_n(S)$ the mass of $\pi_0^*(S)$ on
$\cap_{0\leq j\leq n}\pi_j^{-1}(M''\times N)$.

\begin{lemma} \label{lemma_lov}
Let $\delta$ be a constant strictly larger than  $\widetilde d_q^+$.
Then  there exists a constant $C>0$  such that for all positive closed
current $S$ of bidegree
$(k-q,k-q)$ supported on $M\times N'$
we have $\volume_n(S)\leq C\delta^n\|S\|_D$.
\end{lemma}
\proof
Observe that $f^j$ can be identified with $\pi_j\circ \pi_0^{-1}$. This
allows to
write $\volume_n(S)$ as the following sum of $(n+1)^q$ integrals
$$\volume_n(S)=\Big\langle \pi_0^*(S),\Big(\sum \pi_j^*(\omega_v)\Big)^q
\Big\rangle = \sum_{0\leq n_j\leq n} \int S\wedge (f^{n_1})^*\omega_v\wedge\ldots\wedge
(f^{n_q})^*\omega_v.$$
Lemma \ref{lemma_pre_lov} applied to a constant $\delta'>\widetilde d_q^+$
implies that 
$\volume_n(S)\leq C'n^q {\delta'}^n\|S\|$ for some constant $C'>0$. We
obtain the result by choosing a $\delta'$ smaller than $\delta$.
\endproof

A subset $A$ of $D$ is said to be {\it $(n,\epsilon)$-separated} if
$f^j$ is defined on $A$ with $f^j(A)\subset D'':=M''\times N''$ for $0\leq
j\leq n$ and for every distinct
points $a$, $b$ in $A$ the distance between $f^j(a)$ and $f^j(b)$ is
larger than $\epsilon$ for at least one $j$ with $0\leq j\leq
n$. Define for a subset $X$ of $D$ the {\it topological entropy} of
$f$ on $X$ by
$$h_X(f):=\sup_{\epsilon >0}\limsup_{n\rightarrow\infty} {1\over n} \log
\max\#\big\{A\subset X,\quad A\quad (n,\epsilon)\mbox{-separated}\big\}.$$ 
We have the following version of the Gromov's inequality, see also \cite{Gromov, ds2, deThelin}.
\begin{proposition} \label{prop_entropy_sub}
Let $\delta$ be a constant strictly larger than  $\widetilde d_q^+$ with $q\leq
p-1$.
Let $X$ be a horizontal subvariety of dimension $q$ of $D$. 
Then for every $\epsilon>0$ there is a constant $C_\epsilon>0$ such that
every $(n,\epsilon)$-separated subset in $X$ contains at most
$C_\epsilon\delta^n$ points. In particular, we have  
$h_X(f)\leq \log \widetilde d_q^+$.
\end{proposition}
\proof
We can choose $N'$ such that $X$ is contained in $M\times N'$. We can
also assume that $\epsilon$ is small enough. So,
$X$ defines a horizontal positive closed current $[X]$ of bidegree
$(k-q,k-q)$. Lemma \ref{lemma_lov} applied to $M'$ instead of $M''$,
implies that the volume of $\pi_0^{-1}(X)$
in $\cap_{0\leq j\leq n}\pi_j^{-1}(M'\times N)$ is smaller than
$C\delta^n$ for some constant $C>0$.

Consider an $(n,\epsilon)$-separated subset $A$ of
$X$. For every $a$ in $A$ denote by $B_a$ the ball of center
$(a,f(a),\ldots, f^n(a))$ and of diameter $\epsilon$ in
$D^{n+1}$. Since $A$ is $(n,\epsilon)$-separated, these balls are
disjoint. Since $\epsilon$ is small and the center of $B_a$ is in 
$\cap_{0\leq j\leq n}\pi_j^{-1}(D'')$, these balls are contained in
$\cap_{0\leq j\leq n}\pi_j^{-1}(M'\times N)$. It follows that the
total volume of $B_a\cap \pi_0^{-1}(X)$ is bounded by $C\delta^n$. On
the other hand, an inequality of Lelong \cite{Lelong} says that the
volume of $B_a\cap \pi_0^{-1}(X)$ is bounded from below by a constant
depending only on $\epsilon$. Hence, the number of the balls $B_a$ is
$\lesssim \delta^n$. This implies that $\#A\lesssim \delta^n$ and
completes the proof.
\endproof

Recall that it is proved in \cite{ds6} that $\mu$ is of maximal
entropy $\log d$. This also holds for $f^{-1}$ since the main
dynamical degree of $f^{-1}$ is also equal to $d$. Let
$B_{-n}(x,\epsilon)$ denote the {\it Bowen $(-n,\epsilon)$-ball} with
center $x$, i.e. the set of the points $y$ such that $f^{-j}(y)$ is
defined and $\|f^{-j}(y)-f^{-j}(x)\|\leq \epsilon$ for $0\leq j\leq
n$.  {\it The entropy $h(\mu)$ for $f^{-1}$} can be obtained by the
following Brin-Katok formula \cite{bk}
$$h(\mu):=\sup_{\epsilon>0}\liminf_{n\rightarrow \infty} -{1\over n}
\log\mu(B_{-n}(x,\epsilon))$$
for $\mu$-almost every $x$. So, for every $\theta>0$, there are positive
constants $C$, $\epsilon$ and a Borel set $\Sigma_0$ with $\mu(\Sigma_0)>3/4$ such that
$\mu(B_{-n}(x,6\epsilon))\leq C e^{-n(\log d-\theta)}$ for
$x\in\Sigma_0$ and $n\geq 0$.

\medskip

\noindent
{\bf Proof of Theorem \ref{th_hyperbolic}.}
Assume in order to reach a contradiction that $\mu$
admits at least $k-p+1$ Lyapounov exponents strictly smaller than
${1\over 2k}\log(d/\widetilde\delta_+)$. 
Let $q\leq p-1$ be an integer and $\lambda < {1\over
  2k}\log(d/\widetilde\delta_+)$ a positive constant such that $\mu$
admits exactly $k-q$ Lyapounov exponents strictly smaller than
$\lambda$ and the other ones are larger than or equal to ${1\over
  2k}\log(d/\widetilde\delta_+)$. 
We are going to construct a complex subspace $F$ of dimension $q$,
contradicting the estimate in Proposition \ref{prop_entropy_sub}, i.e. with too many
$(n,\epsilon)$-separated points.

Fix a positive constant $\theta$ such that $\theta\ll \lambda$ and $\theta\ll  {1\over
  2k}\log(d/\widetilde\delta_+) -\lambda$. 
By Oseledec-Pesin theory (replacing $f$ by an iterate $f^n$ and
$\theta$, $\lambda$, $d$, $\widetilde\delta_+$ by $n\theta$,
$n\lambda$, $d^n$, $\widetilde\delta_+^n$ if necessary), 
we can assume that there is  a decomposition
$T_x=E_x\oplus F_x$ 
for $\mu$ almost every $x$ with the following properties: 
\begin{enumerate}
\item[-] $E_x$ and $F_x$ are vector spaces of dimension $k-q$ and $q$
  respectively.
\item[-] The vector bundles $E_x$ and $F_x$ are $f$-invariant.
\item[-] There is a Borel set $\Sigma\subset \Kc$ with
  $\mu(\Sigma)\geq 1/2$ and a constant $\eta>0$ such
      that 
$$\|Df^{-1}(v)\|\geq  e^{-\lambda}\|v\|,\ \ 
  \|Df^{-1}(u)\|\leq  e^{-\lambda-7\theta}\|u\|,\ \
\measuredangle\big(E_{f^{-n}(x)},F_{f^{-n}(x)}\big)\geq \eta e^{-n\theta}$$ 
for $v\in E_x$, $u\in F_x$,  $x\in \Sigma$ and $n\geq 0$.
\end{enumerate}

We now identify each $T_x$ with $\C^k$ and consider $x$ as the
origin. Fix coordinate systems on $E_x$ and $F_x$ so that the
associated distances coincide with the distances induced by the
standard metric on $\C^k$. On $T_x=E_x\oplus F_x$ we use the
coordinate system induced by the fixed coordinates on $E_x$ and
$F_x$. We call it {\it dynamical coordinate system}.
Note that the angle between $E_x$ and $F_x$, with respect
to the standard coordinates, might be
small and in this case there is a big distorsion of the dynamical
coordinates with respect to the standard ones. 

Fix a positive constant $c$ small enough, $c\ll\eta$ and $c\ll\epsilon$ where 
$\epsilon$ is the constant associated to $\theta$ as above. 
Let $B_{x_{-n}}$ denote the (small)
ball of
radius $ce^{-n(\lambda+7\theta)}$ of center $x_{-n}:=f^{-n}(x)$ in
$E_{x_{-n}}$. 
We are interested in graphs in  $T_{x_{-n}}=E_{x_{-n}}\oplus F_{x_{-n}}$ of
holomorphic maps over $B_{x_{-n}}$.

\medskip

\noindent
{\bf Claim 1.} {\it For every $x\in
  \Sigma$ there are holomorphic maps $h_n:B_{x_{-n}}\rightarrow
  F_{x_{-n}}$ with graph $V_{x_{-n}}$ such that 
 $h_n(0)=0$, $\|Dh_n\|\leq e^{-4n\theta}$ and $f$ sends $V_{x_{-n-1}}$
 into $V_{x_{-n}}$. 
}

\medskip

The proof of this claim is by induction. For $n=0$, it is enough to
choose $h_0=0$. We will obtain $V_{x_{-n}}$ as an open set in
$f^{-1}(V_{x_{-n+1}})$. Consider the map $f^{-1}$ on a small
neighbourhood of $x_{-n+1}$ with image in a neighbourhood of
$x_{-n}$.  In dynamical coordinates for $T_{x_{-n+1}}$ and $T_{x_{-n}}$
we can write
$$f^{-1}(z) = l(z)+r(z)\quad \mbox{with}\quad l=(l',l'') \quad
\mbox{and}\quad r=(r',r'')$$
where $l(z)$ is the linear part of $f$, i.e. the differential $Df^{-1}$ at
$x_{-n+1}$, and $r(z)$ is the rest which is of order $\geq 2$ with
respect to $z$. 

We have $l':E_{x_{-n+1}}\rightarrow E_{x_{-n}}$ and   $l'':F_{x_{-n+1}}
\rightarrow F_{x_{-n}}$. We also have $\|l'(z')\|\geq e^{-\lambda}\|z'\|$ for
$z'\in E_{x_{-n+1}}$ and  
$\|l''(z'')\|\leq e^{-(\lambda+7\theta)}\|z''\|$ for  $z''\in F_{x_{-n+1}}$.
In the standard coordinates, 
the derivatives of $f^{-1}$ are
bounded. Taking into account the distortions
of dynamical coordinates, we have $\|Dr(z)\|\leq
Ae^{6n\theta}\|z\|$ with $A>0$ independent of $c, n, \theta$. Now,
consider two points $z=(z',z'')$ and $w=(w',w'')$ in $E_{x_{-n+1}}\oplus
F_{x_{-n+1}}$ which are contained in $V_{x_{-n+1}}$.  So, 
$\|z\|$ and $\|w\|$ are smaller than $2ce^{-(n-1)(\lambda+7\theta)}$.
Write  $\widetilde z:=(\widetilde z',\widetilde z'')=f^{-1}(z)$ and 
$\widetilde w:=(\widetilde w',\widetilde w'')=f^{-1}(w)$. We deduce
from the estimates on $l'$, $Dr$ and $Dh_{n-1}$ that
\begin{eqnarray*}
\|\widetilde z'-\widetilde w'\| & \geq &
\|l'(z')-l'(w')\|-\|r'(z)-r'(w)\| \\
& \geq & e^{-\lambda}\|z'-w'\| -
2Ae^{6n\theta}ce^{-(n-1)(\lambda+7\theta)}\|z-w\| \\
& \geq & e^{-\lambda}\|z'-w'\| -
4Ae^{6n\theta}ce^{-(n-1)(\lambda+7\theta)}\|z'-w'\|.
\end{eqnarray*}
Hence, $\|\widetilde z'-\widetilde w'\|\geq e^{-(\lambda+\theta)}\|z'-w'\|$
since $c$, $\theta$ are small and $\theta\ll \lambda$. It follows that
$f^{-1}(V_{x_{-n+1}})$ is a graph of a holomorphic map $h_n$ over an open set $B$ of
$E_{x_{-n}}$. The last estimate for $w'=0$ implies that $B$
contains the ball $B_{x_{-n}}$.

On the other hand, we have
\begin{eqnarray*}
\|\widetilde z''-\widetilde w''\| & \leq &
\|l''(z'')-l''(w'')\|+\|r''(z)-r''(w)\| \\
& \leq & e^{-(\lambda+7\theta)}\|z''-w''\| +
2Ae^{6n\theta}ce^{-(n-1)(\lambda+7\theta)}\|z-w\| \\
& \leq & e^{-(\lambda+7\theta)}e^{-4(n-1)\theta}\|z'-w'\| +
4Ae^{6n\theta}ce^{-(n-1)(\lambda+7\theta)}\|z'-w'\|.
\end{eqnarray*}
Therefore, $\|\widetilde z''-\widetilde w''\|\leq e^{-4n\theta}
\|\widetilde z'-\widetilde w'\|$ since $\theta\ll \lambda$ and $c$ is small.
It follows that $\|Dh_n\|\leq e^{-4n\theta}$ and this finishes the proof
of the claim.

\medskip

Note that all the constructed graphs are small and contained in a
small neighbourhood $\Uc$ of  the filled
Julia set $\Kc$. We now come back to the standard metric on
$\C^k$. Let $F_x'$ denote the orthogonal of $E_x$. We use 
coordinate systems on $F_x'$ which induce the standard metric. 
Let $B_{x_{-n}}'$ denote the ball of center $0$ and of radius
$c'e^{-n(\lambda+10\theta)}$ in $E_{x_{-n}}$ with $c'>0$ small enough.
We claim that $V_{x_{-n}}$ contains some flat graph $V_{x_{-n}}'$.

\medskip

\noindent
{\bf Claim 2.}  {\it For every $x\in
  \Sigma$, $V_{x_{-n}}$ contains the graph $V'_{x_{-n}}$ 
of a holomorphic map $h_n':B'_{x_{-n}}\rightarrow
  F'_{x_{-n}}$ such that 
 $h'_n(0)=0$ and $\|Dh'_n\|\lesssim e^{-n\theta}$.
}

\medskip

With the considered coordinates on $E_{x_{-n}}$, $F_{x_{-n}}$ and
$F'_{x_{-n}}$, denote by 
$\tau:E_{x_{-n}}\oplus F_{x_{-n}}\rightarrow E_{x_{-n}}\oplus
F'_{x_{-n}}$ the linear map of coordinate change. Since 
the angle between $E_{x_{-n}}$ and $F_{x_{-n}}$ is larger than
$\eta e^{-n\theta}$, we can write $\tau=(\tau',\tau'')$ with
$\|\tau'(z)-z'\|\lesssim e^{n\theta}\|z''\|$ and $\|\tau''(z)\|\leq \|z''\|$ for
$z=(z',z'')$ in $E_{x_{-n}}\oplus F_{x_{-n}}$. Claim 2 is proved 
using analogous estimates as in Claim 1 where we replace $f^{-1}$ by
$\tau$. We will not give the details here.

We continue the proof of Theorem \ref{th_hyperbolic}. 
Let $A$ be a subset of  $\Sigma\cap\Sigma_0$ such that the balls
 $B_{-n}(x,3\epsilon)$ with centers $x\in A$ are disjoint.  We choose
 $A$ maximal satisfying this property.
So, the balls $B_{-n}(x,6\epsilon)$ with centers $x\in A$ cover
$\Sigma\cap\Sigma_0$. 
Since $\mu(\Sigma\cap\Sigma_0)\geq 1/4$ and $\mu(B_{-n}(x,6\epsilon))\leq
Ce^{-n(\log d-\theta)}$, $A$ contains at least 
$(4C)^{-1} e^{n(\log d -\theta)}$ points. Consider the graphs
$V_{x_{-n}}$ and $V_{x_{-n}}'$ constructed above for $x\in A$. 
Since the balls $B_{-n}(x,3\epsilon)$ are disjoint, the set of
$x_{-n}$ are $(n,3\epsilon)$-separated. Claim 1 implies that 
the diameter of $V_{x_{-n}}$ is
 smaller than $\epsilon$.
So, if we
replace each $x_{-n}$ by a point $x_{-n}'$ in $V_{x_{-n}}$ the resulting
set is always
$(n,\epsilon)$-separated.

Let $\Pi$ be an orthogonal projection of $\C^k=\C^p\times \C^{k-p}$ onto a subspace $E$ of
dimension $k-q$. If $E$ is a product of 
a subspace of
$\C^p$ with $\C^{k-p}$, then the fibers of $\Pi$ which are close enough to $\Kc$
(in particular the fibers which intersect $\Uc$) are horizontal in
$D$. This property holds for the projection on any small perturbation
of $E$. So, we can choose a finite number of projections
$\Pi_1$, $\ldots$, $\Pi_N$ on $E_1$, $\ldots$, $E_N$
satisfying this property,  and a constant $\theta_0>0$ such that any
subspace $F$ of dimension $q$ in $\C^k$ has an angle $\geq \theta_0$
with at least one of $E_i$. We deduce from Claim 2 that for each
of the considered graphs $V'_{x_{-n}}$, the volume of 
$\Pi_i(V'_{x_{-n}})$ is $\geq c'' e^{-2n(k-q)(\lambda+10\theta)}$ for at
  least one projection $\Pi_i$ with a fixed constant $c''>0$.  Choose
  an $i$ such that this property holds for at least $N^{-1}\#A$ 
  graphs $V'_{x_{-n}}$. Since $\#A\geq (4C)^{-1} e^{n(\log d
    -\theta)}$, the sum of the volumes of $\Pi_i(V_{x_{-n}})$ is
  $\gtrsim e^{n(\log d
    -\theta)-2n(k-q)(\lambda+10\theta)}$.
Hence, there is a fiber $F$ of $\Pi_i$ which intersects $\gtrsim e^{n(\log d
    -\theta)-2n(k-q)(\lambda+10\theta)}$ graphs $V_{x_{-n}}$. It
  follows that $F$ contains an $(n,\epsilon)$-separated subset of 
$\gtrsim e^{n(\log d
    -\theta)-2n(k-q)(\lambda+10\theta)}\geq e^{n(\log \widetilde
    \delta_+ + \theta)}$ points since $\theta\ll {1\over
  2k}\log(d/\widetilde\delta_+) -\lambda$. This contradicts
  Proposition \ref{prop_entropy_sub} for $X=F$ since $\widetilde \delta_+\geq
  \widetilde d_q^+$,
and finishes the proof of Theorem \ref{th_hyperbolic}.
\hfill $\square$

\begin{remark} \rm
The above bound $ {1\over  2k}\log(d/\widetilde\delta_+)$ can be
replaced by the infimum of the numbers
${1\over  2(k-q)}\log(d/\widetilde d_q^+)$
for $q\leq p-1$. 
\end{remark}

\begin{remark}\rm
The fact that we are in the holomorphic setting is used only in
Proposition \ref{prop_entropy_sub} in order to get an estimate on the
topological entropy on
analytic manifolds of dimension $q$. The result still holds for real
$\Cc^{1+\alpha}$ horizontal-like maps (i.e. non-uniformly hyperbolic
horshoes) with an ergodic invariant measure with compact
support. We only need that the entropy of the measure is strictly
larger than the entropy on vertical subspaces of dimension $\leq
k-p-1$ and horizontal manifolds of dimension $\leq p-1$, see also
Newhouse, Buzzi and de Th\'elin \cite{ne, deThelin}.
\end{remark}


\section{Examples and open problems} \label{section_ex}

Consider a polynomial automorphism $f$ of $\C^k$. We extend $f$ to
a birational map on the projective space $\P^k$. Let $I_+$ and $I_-$
denote the indeterminacy sets of $f$ and $f^{-1}$. They are in the
hyperplane at infinity $L_\infty:=\P^k\setminus \C^k$ and we assume that they are non-empty. 
When $I_+$ and $I_-$ have empty intersection,
$f$ is said to be {\it regular}. This class of automorphisms was
introduced and studied in \cite{si}. In dimension $k=2$, they are 
the H{\'e}non type maps and any polynomial automorphism  of positive
entropy is conjugated to a regular automorphism.

There is an integer $p$ such that $\dim I_+=k-p-1$ and $\dim
I_-=p-1$. If $d_+$ and $d_-$ denote the algebraic degrees of $f$ and
$f^{-1}$, we have $d_+^p=d_-^{k-p}$. 
At infinity we have $f(L_\infty\setminus
I_+)=I_-$ and $f^{-1}(L_\infty\setminus I_-)=I_+$.
Define the {\it filled Julia  sets} by
$$\Kc_+:=\big\{z\in\C^k,\ (f^n(z))_{n\geq 0} \mbox{ bounded in } \C^k\big\}$$
and
 $$\Kc_-:=\big\{z\in\C^k,\ (f^{-n}(z))_{n\geq 0} \mbox{ bounded in } \C^k\big\}.$$
These sets are invariant under $f^{-1}$, $f$ and satisfy
$\overline\Kc_+=\Kc_+\cup I_+$,  $\overline\Kc_-=\Kc_-\cup I_-$.
One associates to $f$ and $f^{-1}$ the following
functions, called {\it Green functions}
$$G^+(z):=\lim_{n\rightarrow \infty} d_+^{-n} \log^+\|f^n(z)\|\qquad
\mbox{and} \qquad G^-(z):=\lim_{n\rightarrow \infty} d_-^{-n} \log^+\|f^{-n}(z)\|,$$
where $\log^+:=\max(\log,0)$. These functions are continuous
p.s.h. on $\C^k$. It follows from \cite[Proposition 2.4]{ds4} that $G^+$
and $G^-$ are H\"older continuous.
They satisfy $G^+\circ f=d_+G^+$ and
$G^-\circ f^{-1}=d_-G^-$. It is shown in \cite{ds8} that the {\it Green currents}
$$T_+:=(\ddc G^+)^p\qquad \mbox{and}\qquad T_-:=(\ddc G^-)^{k-p}$$
are, up to a multiplicative constant, the unique positive closed
currents of bidegrees $(p,p)$ and $(k-p,k-p)$ with
support in $\Kc_+$ and $\Kc_-$ respectively. These currents are
invariant: $f^*(T_+)=d_+^pT_+$ and $f_*(T_-)=d_-^{k-p}T_-$. Note that
to prove the uniqueness we do not assume invariance.

The family of regular automorphisms is large but for simplicity we
restrict to the case where 
the indeterminacy sets $I_+$ and $I_-$ are linear. In what follows, we
assume that 
$$I_+=\{z_0=z_1=\cdots=z_p=0\}\qquad \mbox{and}\qquad 
I_-=\{z_0=z_{p+1}=\cdots=z_k=0\}$$ 
where $[z_0:\cdots:z_k]$ denotes
  the homogeneous coordinates of $\P^k$, $\C^k$ is identified to the
  chart $\{z_0=1\}$ and the hyperplane at infinity $L_\infty$ is given by the
  equation $z_0=0$. The following proposition allows to apply the
  results in the previous sections to the small (possibly transcendental)
  pertubations of $f$ and proves Corollary \ref{cor_ex}.

\begin{proposition} \label{prop_ex}
Let $f$ be a regular polynomial automorphism of $\C^k$ as above. 
Let $B_s^R$ denote the ball of center $0$ and of radius $R$ in $\C^s$.
Then, if $R$ is large enough, 
any holomorphic map $f_\epsilon$ on $B_p^R\times B_{k-p}^R$, close enough to 
$f$, is horizontal-like with the main dynamical degree
$d=d_+^p=d_-^{k-p}$. Moreover, $d$ is strictly
larger than the other dynamical degrees associated to $f_\epsilon$ and $f_\epsilon^{-1}$.
\end{proposition}
\proof
By Proposition \ref{prop_pertub}, it is enough to check that $f$ restricted to
 $B_p^R\times B_{k-p}^R$ is a horizontal-like map of main dynamical degree $d$ which
is strictly
larger than the other dynamical degrees. Write, using the coordinates
$(z_1,\ldots,z_k)$ of $\C^k$ 
$$f=(f',f'') \quad \mbox{with}\quad f'=(f_1,\ldots,f_p)\quad
\mbox{and}\quad f''=(f_{p+1},\ldots,f_k).$$
Since $f(L_\infty\setminus I_+)=I_-$, 
the equation of $I_-$ implies that the components of $f''$ have degree
$\leq d_+-1$ and the components of $f'$ have degree $d_+$. Moreover,
if $f_j^+$ denotes the homogeneous part of degree $d_+$ of $f_j$, the
equation of $I_+$ implies that $f_1^+=\cdots=f_p^+=0$ only when
$z_1=\cdots=z_p=0$. 
The restriction of $f$ to $I_-$ defines an endomorphism of algebraic
degree $d_+$.

Since
$R$ is large, it follows that $\|f'(z)\|>R$ for $z$ in the vertical
boundary of $B_p^R\times B_{k-p}^R$. Hence, $f^{-1}(B_p^R\times B_{k-p}^R)$ does
not intersect the vertical boundary of $B_p^R\times B_{k-p}^R$. In the
same way, we show that  $f(B_p^R\times B_{k-p}^R)$ does
not intersect the horizontal boundary of $B_p^R\times B_{k-p}^R$. This
proves that $f$ restricted to $B_p^R\times B_{k-p}^R$ is
horizontal-like. In order to avoid confusion, let us denote by
$\overline f$ the horizontal-like map on $D:=B_p^R\times B_{k-p}^R$ associated to $f$.

Since $\overline \Kc_+=\Kc_+\cup I_+$, the
equation of $I_+$ implies that  $\Kc_+$ restricted to $D$ 
is vertical. The restriction of $T_+$ to $D$
is vertical and invariant under $d^{-1}\overline f^*$. So, the main dynamical degree
of $\overline f$ is equal to $d$. It remains to check that the other dynamical
degrees are strictly smaller than $d$.

Fix an $\alpha>0$ small enough so that $\overline f^{-1}(D)\subset
B_p^{R-2\alpha}\times B_{k-p}^R$ and $\overline f(D)\subset
B_p^R\times B_{k-p}^{R-2\alpha}$. So $\overline f$ is
horizontal-like on $D':=B_p^{R-\alpha}\times B_{k-p}^{R-\alpha}$
and on $D'':=B_p^{R-2\alpha}\times B_{k-p}^{R-2\alpha}$.
Consider the family $\Qc_h$ of horizontal positive closed currents of bidegree
$(k-s,k-s)$ and of mass 1 in $D''$ with $s\leq p-1$. 
We will show
that the mass of $(\overline f^n)_*S$ on
$D''$  for $S\in\Qc_h$, is of order $O(d_+^s)$. 
This implies that the dynamical degree
$d^+_s$ of $\overline f$ is $\leq d_+^s$ and then is strictly smaller
than $d$. The proof is analogous for the degrees $d_s^-$ associated to
$\overline f^{-1}$.

Observe that $S':=\overline f_*(S)$ is horizontal in $D'$ and
has bounded mass.
Let $\omega_\FS:=\ddc H$, with $H:=\log(1+\|z\|^2)^{1/2}$, be the
Fubini-Study form on $\P^k$.  Since the standard K\"ahler form on $\C^k$ and $\omega_\FS$ are
comparable in compact sets of $\C^k$, it is enough to estimate the
mass of $\omega_\FS^s \wedge (\overline f^n)_*S$ on
$D''$. We have 
\begin{equation} \label{eq_mass_FS}
\int_{D''} 
\omega_\FS^s \wedge (\overline f^n)_*S = \int_{\overline
  f^{-n+1}(D'')} (f^{n-1})^*\omega_\FS^s\wedge S'\leq
\int_{D'} (f^{n-1})^*\omega_\FS^s\wedge S'
\end{equation}
since $\overline f^{-n+1}(D'')\subset
B_p^{R-2\alpha}\times B_{k-p}^R$ and $\supp(S')\subset B_p^R\times
B_{k-p}^{R-\alpha}$. 
It was shown in \cite{si} that $d_+^{-n} \log^+\|f^n(z)\|$ converge
locally 
uniformly to $G^+$. We deduce easily that $d_+^{-n}H\circ f^n$
converge also locally uniformly to $G^+$. It follows from the theory of
intersection of currents, see \cite{Demailly, fs2} that the family of
currents 
$$d_+^{-sn} (f^n)^*\omega_\FS^s\wedge S = d_+^{-sn} (\ddc H\circ
f^n)^s\wedge S$$
is relatively compact. Hence, the integrals in (\ref{eq_mass_FS}) are
$\lesssim d_+^{sn}$ and 
the mass of $(\overline
f^n)_*S$ on $D''$ is $\lesssim d_+^{sn}$. This completes the proof.
\endproof

\begin{remark}\rm
The restriction of $\Kc_+$, $\Kc_-$, $T_+$ and $T_-$ to $D=B_p^R\times
B_{k-p}^R$ coincide with the filled Julia sets and the Green currents
constructed for $\overline f$. Note that in the context of
horizontal-like maps, $T_+$ is not the unique positive closed $(p,p)$-current with support in
$\Kc_+$. For the horseshoes, this current can be decomposed into
currents of integration on vertical submanifolds of $D$.   
\end{remark}

Many questions have to be considered in the context of horizontal-like maps even
when we assume that the condition on the dynamical degrees is satisfied. We refer to the paper by
Dujardin \cite{duj} for the case of dimension 2, see also
\cite{BedfordLyubichSmillie, dds}.

\begin{question}\rm
Let $f$ be an invertible horizontal-like map as above. Is the sequence
$(d_s^+)_{0\leq s\leq p}$ of dynamical degrees of $f$ increasing?
\end{question}

\begin{question} \rm
Let $f$ be an invertible horizontal-like map as above. Is the
Green current $T_+$ laminar? More precisely, is it decomposable into
currents of integration on complex manifolds, not necessarily closed,
in $D$?
\end{question} 

We refer to \cite{Dinh1, deThelin1} for recent results on laminar
currents in higher dimension. The following
problems are also open for regular polynomial automorphisms.

\begin{question} \rm
Is the
equilibrium measure $\mu$ the intersection in the geometrical sense of $T_+$ and
$T_-$? More precisely, is it possible to decompose $T_+$ and $T_-$
into currents of integration on complex manifolds and to obtain $\mu$ as
an average on the intersections of such manifolds?
\end{question}

\begin{question} \rm
Are saddle periodic points
equidistributed with respect to $\mu$? It is not difficult to show
that there are $d^n$ periodic points of period $n$ counted with multiplicities. 
\end{question}

\begin{question} \rm
Is the Hausdorff dimension of $\mu$ positive? Is there a
relation between this dimension and the Lyapounov exponents of $\mu$?
\end{question}

In the case of regular polynomial automorphisms, since $\mu=(\ddc
G^+)^p\wedge (\ddc G^-)^{k-p}$ and $G^+$, $G^-$ are H\"older
continuous, $\mu$ gives no mass to sets of small Hausdorff dimension,
see e.g. \cite[Th\'eor\`eme 1.7.3]{si}.

We refer to Dupont \cite{Dupont}, Ledrappier-Young \cite{LedrappierYoung} 
and the references therein for analogous problems in other contexts.

The dependence of Lyapounov exponents on the map can be studied
following the works by Bassanelli-Berteloot \cite{BassanelliBerteloot} and
Pham \cite{Pham}.


\noindent
T.-C. Dinh, UPMC Univ Paris 06, UMR 7586, Institut de
Math{\'e}matiques de Jussieu, F-75005 Paris, France. {\tt
  dinh@math.jussieu.fr}, {\tt http://www.math.jussieu.fr/$\sim$dinh} 

\

\noindent
V.-A.  Nguy{\^e}n,
Mathematics Section,
the Abdus Salam International Centre 
 for Theoretical Physics,
Strada costiera, 11, 34014 Trieste, Italy.
{\tt vnguyen0@ictp.trieste.it}

\

\noindent
N. Sibony, Math{\'e}matique-B{\^a}timent 425, UMR 8628, Universit{\'e} Paris-Sud,
91405 Orsay, France. {\tt Nessim.Sibony@math.u-psud.fr}


\small

\begin{thebibliography}{99}

\bibitem{BassanelliBerteloot}
G. Bassanelli, F. Berteloot, Bifurcation currents in holomorphic
dynamics on ${\bf P}^k$, {\it  J. Reine Angew. Math.}, {\bf  608}  (2007), 201-235. 


\bibitem{BedfordLyubichSmillie}
E. Bedford, M. Lyubich, J. Smillie,
Polynomial diffeomorphisms of $\C^2$
  V: The measure of maximal entropy and laminar currents,
  \textit{Invent. Math.}, \textbf{112(1)} (1993), 77-125.

\bibitem{bc}
M. Benedicks, L. Carleson, The dynamics of the H{\'e}non map,  
{\it  Ann. of Math. (2),}  {\bf 133}  (1991),  no. 1, 73-169. 

\bibitem{BonattiDiazViana}
C. Bonatti, L.J. D{\'\i}az, M. Viana,  {\it Dynamics beyond uniform
  hyperbolicity. A global geometric and probabilistic perspective,} 
Encyclopaedia of Mathematical Sciences, {\bf 102}, Mathematical Physics, III. Springer-Verlag, Berlin, 2005. 


\bibitem{bt}
R. Bott, L.W. Tu, {\it Differential forms in algebraic topology,}
 Graduate Texts in Mathematics, {\bf 82}, Springer-Verlag, New York-Berlin, 1982.


\bibitem{bk}
M. Brin, A. Katok,  On local entropy, 
  Geometric dynamics (Rio de Janeiro, 1981),  30-38, {\it Lecture Notes in
  Math.}, {\bf 1007}, Springer, Berlin, 1983. 

\bibitem{Demailly}
J.-P. Demailly, {\it Complex analytic geometry}, available at \\
{\tt http://www.fourier.ujf-grenoble.fr/$\sim$demailly} 

 \bibitem{deThelin1}  H. de Th{\'e}lin, Un crit{\`e}re de laminarit{\'e}
   locale en dimension quelconque, {\it Amer. J. Math.},  {\bf 130}  (2008),  no. 1, 187-205.

 \bibitem{deThelin}  -----,  Sur les exposants de Lyapounov des applications m{\'e}romorphes,
   {\it Invent. math.}, {\bf 172} (2008), no. 1, 89-116.

\bibitem{Dinh1}
T.-C. Dinh, Suites d'applications m{\'e}romorphes multivalu{\'e}es et
courants laminaires, {\it J. Geometric Analysis},  {\bf 15} (2005), no. 2, 207-227.

\bibitem{din} 
-----,  Decay of correlations for H{\'e}non maps,
      {\it Acta Math.},  {\bf  195}  (2005), 253-264.

\bibitem{dds}
T.-C. Dinh, R. Dujardin, N. Sibony, 
 On the dynamics near infinity of some
polynomial mappings in $\C^2$, {\it Math. Ann.}, {\bf 333}, No. 4
(2005), 703-739.

\bibitem{ds1}
T.-C. Dinh, N. Sibony, 
Dynamique des applications d'allure
polynomiale, {\it J. Math. Pures  Appl.}, {\bf 82} (2003),
367-423. 


\bibitem{ds2}
-----, Une borne sup{\'e}rieure pour
l'entropie topologique  d'une application rationnelle, {\it Ann. of
  Math.},  {\bf 161}  (2005),  no. 3, 1637-1644. 

\bibitem{ds4}
-----, Green currents for holomorphic automorphisms of 
compact K{\"a}hler manifolds, \textit{J. Amer. Math. Soc.}, {\bf 18} (2005), 291-312.

\bibitem{ds5}
-----, Dynamics of regular birational maps in $\mathbb{P}^k$, {\it J. Funct. Anal.}, 
{\bf 222} (2005), no 1, 202-216.

\bibitem{ds6}  
-----,
  Geometry of currents, intersection theory and dynamics of horizontal-like maps,
  {\it Ann. Inst. Fourier (Grenoble)},  {\bf  56}  (2006),  no. 2, 423--457. 
  
  \bibitem{ds7} 
-----,
Pull-back of currents by holomorphic maps,  
{\it Manuscripta Math.}, {\bf  123}  (2007),  no. 3, 357-371.


\bibitem{ds8}  -----, 
Super-potentials of positive closed currents, intersection theory and
dynamics, {\it Acta Math.}, to appear. {\tt arXiv:math/0703702}


\bibitem{duj} 
R. Dujardin, H{\'e}non-like mappings in $\C^2$, {\it Amer. J. Math.}, 
{\bf 126} (2004), 439-472.

\bibitem{Dupont}
C. Dupont, A lower bound for the dimension of the maximal entropy
measure of endomorphisms of $\C\P^2$, {\it preprint}, 2007.

\bibitem{fs1}
J.-E. Forn\ae ss, N. Sibony, Complex H{\'e}non mappings in $\C^2$
and Fatou-Bieberbach domains, \textit{Duke Math. J.}, \textbf{65}
(1992), 345-380.

\bibitem{fs2}
-----, Oka's inequality for currents and
applications, \textit{Math. Ann.}, \textbf{301} (1995), 399-419.

\bibitem{Gromov}
M. Gromov,  On the entropy of holomorphic maps,
\textit{Enseignement Math.}, \textbf{49} (2003), 217-235. {\it Manuscript} (1977).

\bibitem{hk1}
G. Henkin, J. Leiterer,  {\it Theory of functions on complex manifolds},
 Monographs in Mathematics, {\bf 79}, Birkh{\"a}user Verlag, Basel, 1984.

\bibitem{hk2}
-----,  {\it Andreotti-Grauert theory by integral formulas},
 Progress in Mathematics, {\bf 74}, Birkh{\"a}user Boston, Inc., Boston, MA, 1988. 


\bibitem{Hormander}
L. H{\"o}rmander, {\it An introduction to complex analysis in several
  variables}, Third edition, North-Holland Mathematical Library, {\bf
  7}, North-Holland Publishing Co., Amsterdam, 1990.


\bibitem{kh}
A. Katok, B. Hasselblatt, {\it Introduction to the modern
theory of dynamical systems}, Cambridge Univ.,
  Press. Encycl. of Math. and its Appl. \textbf{54}, 1995. 

\bibitem{LedrappierYoung}
F. Ledrappier, L.-S. Young,
The metric entropy of diffeomorphisms. II. Relations between entropy, exponents and dimension,
{\it Ann. of Math. (2)}, {\bf 122} (1985), no. 3, 540-574. 


\bibitem{Lelong}
P. Lelong, {\it Fonctions plurisousharmoniques et formes
diff{\'e}rentielles positives}, Dunod Paris, 1968.


\bibitem{ne}
S.E. Newhouse,  Entropy and volume,  {\it  Ergodic Theory
  Dynam. Systems,} {\bf 8}  (1988),  
Charles Conley Memorial Issue, 283-299.

\bibitem{pe}
Y.B. Pesin, Characteristic Liapunov exponents, and smooth ergodic theory, 
{\it Russian Math. Surveys,} {\bf  32} (1977), no. 4, 55-114.

\bibitem{Pham}
N.-M. Pham, Lyapunov exponents and bifurcation current for
polynomial-like maps, {\it preprint}, 2005. {\tt arXiv:math/0512557}


\bibitem{ru}
W. Rudin, {\it Function theory in the unit ball of $\C^{n}$},
 Grundlehren der Mathematischen Wissenschaften, \textbf{241},
 Springer-Verlag, New York-Berlin, 1980. 

\bibitem{si}
N. Sibony, Dynamique des applications rationnelles de
$\mathbb{P}^k$, in  {\it Dynamique et g{\'e}om{\'e}trie  complexes}  
(Lyon, 1997), 97-185, \textit{Panoramas et Synth{\`e}ses},  \textbf{8},
Soc. Math. France, Paris, (1999).

\bibitem{tri}  H.  Triebel, {\it Interpolation Theory, Functions Spaces, Differential Operators,}
North-Holland Math. Library, \textbf{18} North-Holland, Amsterdam--New  York, 1978.


\bibitem{Vigny}
G. Vigny, Dynamics semi-conjugated to a subshift for some polynomial
mappings in $\Bbb C\sp 2$,  {\it Publ. Mat.},  {\bf 51}  (2007),  no. 1, 201-222. 


\bibitem{yo} L.-S. Young,  Statistical properties of dynamical systems with some hyperbolicity,
{\it  Ann. of Math. (2),}  {\bf 147}  (1998),  no. 3, 585-650. 

\bibitem{yo2}   -----,  Ergodic theory of chaotic dynamical systems,
{\it  XIIth International Congress of Mathematical Physics (ICMP '97)
  (Brisbane),}  131-143, Int. Press, Cambridge, MA, 1999. 
\end{thebibliography}
\end{document}